\documentclass[11pt, reqno, a4paper]{amsart}
\usepackage{amsmath,amsfonts,amssymb,amscd,amsthm}
\usepackage{hyperref}
\usepackage{bm}
\usepackage{braket}
\usepackage{mathrsfs,centernot,stmaryrd,color,tensor,wasysym, bbm}
\usepackage[truedimen,margin=25truemm]{geometry}
\usepackage{tcolorbox}
\usepackage{mathtools}
\usepackage{tikz}
\usepackage{ulem}
\usepackage{cancel}
\usepackage[sort,compress]{cite} 
\usetikzlibrary{intersections,calc,arrows.meta}
\mathtoolsset{showonlyrefs}

\allowdisplaybreaks

\newcommand{\N}{\mathbb{N}}
\newcommand{\Z}{\mathbb{Z}}

\newcommand{\R}{\mathbb{R}}

\newcommand{\abs}[1]{\left| {#1} \right|}

\newcommand{\nor}[1]{\left\| {#1} \right\|}

\theoremstyle{plain}
\newtheorem{th.}{Theorem}[section]
\newtheorem{prop.}[th.]{Proposition}
\newtheorem{lem.}[th.]{Lemma}
\newtheorem{cor.}[th.]{Corollary}

\theoremstyle{definition}
\newtheorem{def.}[th.]{Definition}
\newtheorem{rmk.}[th.]{Remark}
\newtheorem{conj.}[th.]{Conjecture}

\newcommand{\BC}{\mathbb{C}}
\newcommand{\BD}{\mathbb{D}}
\newcommand{\BE}{\mathbb{E}}

\newcommand{\BN}{\mathbb{N}}

\newcommand{\BP}{\mathbb{P}}

\newcommand{\BR}{\mathbb{R}}

\newcommand{\BT}{\mathbb{T}}

\newcommand{\CB}{\mathcal{B}}

\renewcommand{\CD}{\mathcal{D}}

\newcommand{\CF}{\mathcal{F}}

\newcommand{\CH}{\mathcal{H}}

\newcommand{\CL}{\mathcal{L}}

\newcommand{\CR}{\mathcal{R}}

\newcommand{\CV}{\mathcal{V}}

\newcommand{\FS}{\mathfrak{S}}

\newcommand{\al}{\alpha}
\newcommand{\be}{\beta}
\newcommand{\ga}{\gamma}
\newcommand{\de}{\delta}
\newcommand{\ep}{\varepsilon}

\newcommand{\ph}{\varphi}
\newcommand{\ta}{\tau}

\newcommand{\si}{\sigma}
\newcommand{\om}{\omega}

\newcommand{\rh}{\rho}
\newcommand{\ch}{\chi}

\newcommand{\De}{\Delta}
\newcommand{\Ph}{\Phi}

\newcommand{\Om}{\Omega}

\newcommand{\pl}{\partial}
\newcommand{\wt}{\widetilde}
\newcommand{\wh}{\widehat}

\newcommand{\supp}{\mathrm{supp}\,}

\newcommand{\loc}{{\rm loc}}

\newcommand{\Dk}[1]{\left[#1\right]}
\newcommand{\K}[1]{\left(#1\right)}

\newcommand{\No}[1]{\left\| #1 \right\|}

\newcommand{\I}{\infty}

\newcommand{\Tr}{{\rm Tr}}

\newcommand{\sd}{\langle \nabla \rangle}

\newcommand{\lxr}{\langle \xi \rangle}

\newcommand{\tx}{{(t,x)}}

\newcommand{\st}{\star}

\newcommand{\ft}{{\frac{1}{4}}}

\newcommand{\tw}{\frac{1}{2}}
\newcommand{\na}{\nabla}
\newcommand{\oc}{\ocircle}
\newcommand{\os}{\oast}
\newcommand{\ls}{\lesssim}

\newcommand{\wo}{{\wt{\om}}}
\newcommand{\wq}{{\wh{q}}}
\newcommand{\ov}{\overline}
\newcommand{\rg}{\rangle}
\renewcommand{\lg}{\langle}

\newcommand{\EQ}[1]{\begin{equation}\begin{aligned} #1 \end{aligned}\end{equation}}
\renewcommand{\supp}{\operatorname{supp}}

\numberwithin{equation}{section}

\makeatletter
\newcommand\thankssymb[1]{\textsuperscript{\@fnsymbol{#1}}}
\makeatother

\title[Probabilistic Strichartz estimates and their applications]
{Probabilistic Strichartz estimates in Schatten classes and their applications to the Hartree equation}
\author{Sonae Hadama\thankssymb{1} \and Takuto Yamamoto\thankssymb{2}}
\thanks{
	\thankssymb{1}
	Research Institute for Mathematical Sciences, Kyoto University, Kita-Shirakawa, Sakyo-ku, Kyoto, Japan 606-8502,
	\texttt{hadama@kurims.kyoto-u.ac.jp}.
	\thankssymb{2}
	Independent Scholar, Shigino-nishi, Joto-ku, Osaka, Japan 536-0014, \texttt{takuto.yamamoto.4525+math@gmail.com}}
\date{\today}

\begin{document}

\maketitle

\begin{abstract}
In this paper, we consider the Strichartz estimates for orthonormal systems in the context of randomization. Frank, Lewin, Lieb, and Seiringer first proved the orthonormal Strichartz estimates. After that, many authors have studied this type of inequality. In this paper, we introduce two randomizations of operators and show that they allow us to treat strictly bigger Schatten exponents than the sharp exponents of the deterministic orthonormal Strichartz estimates. In the proofs, the orthogonality does not have any essential role, and randomness works instead of it. We also prove that our randomizations of operators never change the Schatten classes to which they originally belong. 
Moreover, we give some applications of our results to the Hartree and linearized Hartree equations for infinitely many particles. First, we construct local solutions to the Hartree equation with initial data in wide Schatten classes. By only using deterministic orthonormal Strichartz estimates, it is impossible to give any solution in our settings. Next, we consider the scattering problem of the linearized Hartree equation. One of our randomizations allows us to treat wide Schatten classes with some Sobolev regularities, and by the other randomization, we can remove all Sobolev regularities. 
\end{abstract}

\vspace{+3mm}


\tableofcontents
\section{Introduction}\label{sec:intro}
\subsection{Strichartz estimates in Schatten classes}
In this paper, we consider the Strichartz estimates in Schatten classes:
\begin{align} \label{eq:SS}
	\|\rh(e^{it\De}\ga_0 e^{-it\De})\|_{L^p_t(\BR;L^q_x(\BR^d))} \ls \|\ga_0\|_{\FS^\al},
\end{align}
where $p,q,\al \in [1,\I]$ and $\rh_A = \rh(A)$ is the density function of a linear operator $A$; that is, $\rh_A(x)=A(x,x)$ and $A(x,y)$ is the integral kernel of $A$.
We denote Schatten $\al$--class by $\FS^\al$;
namely, for a compact operator $A \in \CB(L^2_x)$, define
\begin{align}
	\|A\|_{\FS^\al}:=
	\begin{cases}
		&\K{\Tr(|A|^{\al})}^{\frac{1}{\al}} \mbox{ if } 1 \le \al <\I,\\
		&\|A\|_{\CB(L^2_x)} \mbox{ if } \al = \I,
	\end{cases}
\end{align}
where $\Tr$ is the trace on $L^2_x$ and $\CB(L^2_x)$ is the space of all bounded linear operators on $L^2_x$.
It is well-known that $\FS^1$ is the trace class, that $\FS^2$ is the Hilbert--Schmidt class, and that the continuous embedding
\begin{align}
	\FS^{\al_1} \subsetneq \FS^{\al_2}
\end{align}
holds for $1 \le \al_1 < \al_2 \le \I$.
See \cite{2005Simon} for more details on Schatten classes. 
If $\ga_0$ is compact and self-adjoint, we have the spectral decomposition
\begin{align}
	\ga_0 = \sum_{n=1}^\I a_n |f_n\rg \lg f_n |,
\end{align}
where $|f\rg \lg g|(\phi) := \lg g | \phi \rg f$.
Therefore, we can rewrite \eqref{eq:SS} as
\begin{align} \label{eq:OS}
	\No{\sum_{n=1}^\I a_n |e^{it\De}f_n|^2}_{L^p_t(\R;L^q_x)} \ls \|a\|_{\ell^\al}.
\end{align}
If $\al=1$, we immediately obtain \eqref{eq:OS} (therefore, also \eqref{eq:SS}) from the triangle inequality and the usual Strichartz estimate:
\begin{align}
	\|e^{it\De}u_0\|_{L^p_t(\R;L^q_x)} \ls \|u_0\|_{L^2_x}.
\end{align}
However, \eqref{eq:OS} holds for some $\al>1$ if $(f_n)_{n=1}^\I$ is an orthonormal system.
Frank, Lewin, Lieb, and Seiringer first proved this result in \cite{2014FLLS}.
After that, Frank and Sabin obtained the optimal range of exponents that \eqref{eq:OS} holds in \cite{2017FS}.
We can summarize these results as
\begin{th.}[\cite{2014FLLS}Theorem 1; \cite{2017FS}Theorem 8] \label{SS}
	Let $d \ge 1$.
	Let $p,q,\al \in [1,\I]$ satisfy
	\begin{align}
		\frac{2}{p}+\frac{d}{q} =d, \quad 1 \le q < \frac{d+1}{d-1}.
	\end{align}
	Then \eqref{eq:SS} and \eqref{eq:OS} hold for $\al = \frac{2q}{q+1}$. 
	Moreover, $\al = \frac{2q}{q+1}$ is sharp in the following sense:
	If $\al > \frac{2q}{q+1}$, then \eqref{eq:SS} and \eqref{eq:OS} fail.
\end{th.}

In the whole space $\BR^d$, many authors have studied this type of Strichartz estimate.
In \cite{2019BHLNS}, Theorem \ref{SS} was extended to initial data with Sobolev regularity.
When $d=1$, it is unknown whether the endpoint estimate \eqref{eq:OS} with $p=2,q=\I,\al <2$ holds or not.
In \cite{2020BLN}, a weaker version of this endpoint estimate
\begin{align}
	\No{\sum_{n=1}^\I a_n |e^{it\pl_x^2}f_n|^2}_{L^{2,\I}_t(\R;L^\I_x(\BR))} \ls \|a\|_{\ell^\al}
\end{align}
was proved for any $\al <2$.
In \cite{2017CHP} and \cite{2018CHP}, the estimate
\begin{align}
	\||\na|^\tw \rh(e^{it\De} \ga_0 e^{-it\De})\|_{L^2_t(\R;H^{s_0}_x)} \ls \|\sd^s \ga_0 \sd^s \|_{\FS^2}
\end{align}
was proved under the appropriate assumptions for the exponents.
In \cite{2020N}, the same type of estimates as \eqref{eq:OS} was studied on torus $\BT^d$.
All results mentioned above are about free Schr\"{o}dinger equation,
but there are similar studies for other equations, for example, wave, Klein--Gordon and fractional Schr\"{o}dinger equations;
see \cite{2017FS} and \cite{2021BLN}.

After the first version of this article was posted on arXiv, there have been significant developments on the orthonormal Strichartz estimates. In \cite{2023H}, the estimate \eqref{eq:SS} was extended to Schr\"{o}dinger operators with \textit{time-dependent} potential. Further generalizations to Schr\"{o}dinger operators with \textit{time-independent} potential can be found in \cite{Hoshiya 2024 JFA, Hoshiya 2025 AHP, Hoshiya 2025 JMP}. See also \cite{Bez et al 2025 TLMS, Senapati et al 2024, Mondal Song 2025, Ghosh Swain 2024}. Finally, we remark that a new kind of bound for the left-hand side of \eqref{eq:SS} was given via the generalized relative entropy in \cite{Hadama Hong 2025}.

In this paper, we are interested in only Schr\"{o}dinger equations. \textit{We improve the Schatten exponents $\al$ in the estimates \eqref{eq:SS} by using randomization of initial operator $\ga_0$.}

\subsection{Hartree equation for infinitely many particles}
In this paper, we apply our randomized estimates to the Cauchy problem of the Hartree equation:
\begin{equation} \tag{NLH} \label{NLH}
\left\{ \,
\begin{aligned}
	&i\pl_t \ga = [-\De + w \ast \rh_\ga, \ga], \quad \ga:\BR \to \CB(L^2_x), \\
	&\ga(0)=\ga_0,
\end{aligned}
\right.
\end{equation}
where $w:\R^d\to\R$ is a given potential, 
$\ast$ is the space convolution in $\R^d$, and $[\cdot,\cdot]$ is the commutator.
As one of the nonlinear Schr\"{o}dinger equations, the Hartree equation for one particle
\begin{align} \label{eq:one}
	i\pl_t u = (-\De+w \ast |u|^2)u, \quad u:\BR \times \BR^d \to \BC
\end{align}
is well-studied.
We can regard \eqref{NLH} as a many-particle version of \eqref{eq:one}.
See \cite[Introduction]{2015LS} for more details of the background of \eqref{NLH}.

We are interested in the initial value problem \eqref{NLH}.
On the one hand, when $\ga_0$ is in the trace class, the well-posedness of this Cauchy problem was already studied in \cite{1974BDF, 1976BDF, 1976C, 1992Z}; and more recently, the small data scattering was proved in \cite{2021PS}.
More precisely, the above results worked in more general settings than this paper's; for example, they include Hartree--Fock or general Kohn--Sham equations. 
On the other hand, it is important to study the case $\ga_0$ is not in the trace class because we know the following proposition:
\begin{prop.} \label{prop:sta}
	Let $d \ge 1$.
	If $w$ is a finite signed Borel measure on $\R^d$ and $f \in L^1_\xi \cap L^\I_\xi$, then
	$\ga_f := f(-i\na)= \CF^{-1} f \CF$ is a stationary solution to \eqref{NLH}.
\end{prop.}  
Note that $\ga_f$ is not even compact unless $f\equiv0$.
Moreover, there are some important $f$ from the physics point of view;
see, for example, (6)--(9) in \cite{2015LS}.

In this paper, we are interested in \eqref{NLH} that initial data $\ga_0$ is given around the stationary solution $\ga_f$.
Let $Q(t)$ be a perturbation from $\ga_f$; that is, $Q(t):= \ga(t)-\ga_f$.
Then we have
\begin{equation} \tag{$f$-NLH} \label{eq:fNLH}
	\left\{ \,
	\begin{aligned}
		&i\pl_t Q = [-\De+w \ast \rh_Q,Q+\ga_f], \quad Q:\BR \to \CB(L^2_x), \\
		&Q(0)=Q_0.
	\end{aligned}
	\right.
\end{equation} 

Now we explain the known result related \eqref{eq:fNLH}.
Lewin and Sabin \textit{first} formulated the Cauchy problem \eqref{eq:fNLH} in \cite{2015LS} and \cite{2014LS}.
In \cite{2015LS}, they proved the local and global well-posedness when $w \in L^1_x \cap L^\I_x$ and $Q_0 \in \FS^\al$ for some $\al > 1$. 
In \cite{2014LS}, it was proved that $Q(t)$ scatters when $\|Q_0\|_{\FS^{4/3}}$ is sufficiently small in two-dimensional space.
Chen, Hong, and Pavlovi\'{c} extended these results.
In \cite{2017CHP}, they proved the global well-posedness when $w$ is the Dirac delta measure and $f(\xi)=\mathbbm{1}_{\{|\xi|^2\le 1\}}$, where $\mathbbm{1}_A$ is the indicator function of $A$. 
This $\mathbbm{1}_{\{|\xi|^2\le 1\}}$ is one of the most important examples in Proposition \ref{prop:sta}.
In \cite{2018CHP}, they showed the small data scattering in $d$--dimensional spaces when $d \ge 3$.
The asymptotic stability result in \cite{2018CHP} was improved in \cite{2023H} when $d=3$. Finally, we note that, in \cite{Nguyen You 2025}, the linearized equation
\begin{equation}
	i\pl_t Q = [-\De, Q] + [w \ast \rh_Q,\ga_f]
\end{equation}
was studied when $w$ is the Coulomb interaction.

One of the most essential tools in the above results is the Strichartz estimates in Schatten classes.
Therefore, the Schatten classes that initial data belong to are determined by the Strichartz estimates.
\textit{In this paper, we give some probabilistic results that treat strictly bigger Schatten exponents than the deterministic Strichartz estimates allow.}
%%%%%
  
\subsection{Randomization for compact operators}\label{sec:randomization}

It is sometimes possible that PDEs are well-posed by suitable randomization even in supercritical case.
Before defining the randomization used in this paper, we recall some previous studies on various types of randomization and PDEs.
See, for example, \cite{BOP2015a,BOP2015b,BOP2019,B94,B96,2019Br,BT08a,BT08b,BTT2013,DLM19,DLM20,M19,NY19,2016OP,Sp22,2009Th,2008Tz,W32,ZF12}.

This probabilistic approach was initiated by Bourgain in \cite{B94,B96}.
Bourgain considered the periodic nonlinear Schr\"odinger equation with randomized initial data in one and two dimensions.
Using the probabilistic setting and deterministic arguments, Bourgain showed that the periodic NLS is well-posed almost surely.
Inspired by these papers, many authors have studied PDEs with randomized initial data or randomized final data.

Burq and Tzvetkov, in \cite{BT08a,BT08b}, treated wave equations with randomized initial values on three-dimensional compact Riemannian manifolds in the supercritical regime.
Considering the randomization using the eigenfunction expansion for $-\Delta$, they obtained some well-posedness results in a probabilistic sense.

Another commonly used type of randomization is Wiener randomization, defined by the unit-scale decomposition of frequency space.
It was named after the Wiener decomposition in \cite{W32} and was introduced by Zhang and Fang in \cite{ZF12}.
Dodson, L\"uhrmann and Mendelson used Wiener randomization for initial data and obtained some probabilistic results for nonlinear Schr\"odinger and wave equations in low regularity spaces in \cite{DLM19,DLM20}.

Some results for randomized final-state problems have been obtained.
In \cite{M19}, Murphy introduced the randomization defined by the unit-scale decomposition of physical space, then proved almost sure existence and the uniqueness of the wave operator for $L^2$--subcritical NLS above the Strauss exponent.
After that, Nakanishi and the second author extended the result below the Strauss exponent and applied it to the NLS system and the Gross-Pitaevskii equation in \cite{NY19}.
Spitz proved the existence and the uniqueness of the wave operator for the Zakharov system in $\R^3$ by the randomization of the physical space and the angular randomization in \cite{Sp22}.

\bigskip

Here, we prepare the randomization used in this paper.
Let $g_n$ be a mean-zero \textit{real-valued} random variable with distribution $\mu_n$ on a probability space $(\Omega,\mathcal{A},\mathbb{P})$ for each $n$.
Moreover, we assume that $(g_n)_{n=1}^\I$ is independent, and that there exists $C>0$ such that for all $\zeta\in\R$ and $n\in\N$,
\EQ{\label{ineq:g-condition}
	\int_{\R}e^{\zeta x}\,d\mu_n(x)\le e^{C\zeta^{2}}
}
holds.
For example, $(g_n)_{n=1}^\I$ is the family of independent mean-zero Gaussian random variables with a bounded variance $\sigma_n>0$
and distribution $\mu_{n}(x)=(2\pi\sigma_{n})^{-1/2}\exp(-x^2/2\sigma_{n})$.
In this case, the left side of \eqref{ineq:g-condition} is equal to $\exp(\sigma_{n}\zeta^2/2)$.
Another example is a random variable with compactly supported distributions.
We need the condition \eqref{ineq:g-condition} for Lemma \ref{lem:LDE}, the large deviation estimate.
Note that
\begin{align} \label{eq:distr g_n}
	\sup_{n\in \BN}\BE[|g_n|^\al] < \I \text{ for any } 1 \le \al < \I.
\end{align}

\begin{def.}[Singular value randomization] \label{def:random-Salpha}
	Let $\ga_0 \in \CB(L^2_x)$ be compact.
	Then we have the singular value decomposition
	\begin{align}
		\ga_0 = \sum_{n=1}^\I a_n |u_n \rangle \langle v_n|,
	\end{align}
	where $a_n \ge 0$, $(u_n)_{n=1}^\I$ and $(v_n)_{n=1}^\I$ are orthonormal systems in $L^2_x$.
	For any $\omega\in\Omega$, we define $\gamma_0^{\omega}$ as
	\EQ{\label{def rand Salpha}
		\gamma_0^\omega:=\sum_{n=1}^\I a_ng_n(\omega)|u_n \rangle \langle v_n|.
	}
	When $\sd^\si \ga_0 \sd^\si$ is compact for $\si \in \BR$, we define
	\begin{align}
		\ga_0^{\om;\si}:= \sd^{-\si}(\sd^\si \ga_0 \sd^\si )^\om \sd^{-\si}.
	\end{align}
\end{def.}

\begin{rmk.}
	If $(g_n)_{n=1}^\I$ does not gather to $0$ as $n \to \infty$, we cannot get any gain of the Schatten integrability by this randomization; that is, $\gamma_0\in\mathfrak{S}^\alpha\setminus\mathfrak{S}^\beta$ implies $\gamma_0^\omega\in\mathfrak{S}^\alpha\setminus\mathfrak{S}^\beta$ for almost all $\om \in \Om$.
	Similarly, we have no gain of regularity.
	See Proposition \ref{prop:appendix lower} for more details.
\end{rmk.}

Let $\ch:=\mathbbm{1}_{[0,1]^d}$. Then, we have
\begin{equation}\label{eq:partition}
	\sum_{k\in\Z^d}\chi_k(\xi)=1,\quad
	\chi_k(\xi):=\chi(\xi-k)
	\end{equation}for almost all $\xi \in \R^d$.
	Define $\Pi_k :=\CF^{-1} \ch_k \CF$ for all $k\in\Z^d$.
	For each $k \in \Z^{d}$, let $h_k$ be a \textit{real-valued} random variable on a probability space $(\wt{\Om},\wt{A},\wt{\BP})$ satisfying the same assumptions for $g_n$.
	Then for any $\wo \in \wt{\Om}$, we define Wiener randomization by
	\begin{align}\label{eq:rand L2}
		f^\wo := \CR^\wo f := \sum_{k\in\Z^d}h_k(\wo)\Pi_k f.
	\end{align}
	Define the \textit{full randomization} of a compact operator $\ga_0=\sum_n a_n |u_n\rg \lg v_n|$ by
	\begin{align}\label{eq:full rand}
		\ga_0^{\om,\wo} := \CR^\wo \ga_0^\om \CR^\wo = \sum_{n=1}^\I a_ng_n(\om) |u_n^\wo \rangle \langle v_n^\wo|.
	\end{align}
	Moreover, if $\sd^\si \ga_0 \sd^\si$ is compact for $\si \in \BR$, then we define
	\begin{align}
		\ga_0^{\om,\wo;\si}:= \sd^{-\si}(\sd^\si \ga_0 \sd^\si )^{\om,\wo} \sd^{-\si}.
\end{align}

\begin{rmk.}
	If we choose the Rademacher sequence as $(h_k)_{k=1}^\I$, we cannot obtain any lower Schatten exponent or higher Sobolev regularity by full randomization. See Proposition \ref{prop:no gain} for more details.
\end{rmk.}

%%%%%

\subsection{Notations}
We collect necessary notations in this paper.
Define $U(t):= e^{it\De}$.
For any Banach space $X$, we write 
\begin{align}
	L^p_T X = L^p_t(0,T;X) := L^p_t([0,T];X), 
	\quad C_T X = C(0,T;X) := C([0,T];X).
\end{align}
For two operators $A$ and $B$, we write $A_\st B := A B A^*$.
We define Schatten--Sobolev norm by
\begin{align}
	\|A\|_{\CH^{s,\al}} := \|\sd^s_\st A\|_{\FS^\al}.
\end{align}
We write $\CH^s := \CH^{s,2}$.

\subsection{Main results}\label{subsec:main theorem}
\subsubsection{Strichartz estimates in Schatten classes}
Define $ABCD$ by the convex hull of $A, B, C, D$ in Figure 1. 
\begin{figure}[htbp]\begin{center}\begin{tikzpicture}[scale=5]
			\def\r{.015}\def\d{3}\def\s{2.2}
			\coordinate[label=left:O](O)at(0,0);
			\coordinate(QL)at(1.1,0);\coordinate(PL)at(0,1.1);
			\draw[<->,>=stealth,semithick](QL)node[right]{$\frac{1}{q}$}--(O)--(PL)node[left]{$\frac{1}{p}$};
			
			\coordinate[label=left:A](A)at(0,{(\d-\s)/2});
			\coordinate[label=below:B](B)at({(\d-\s)/\d},0);
			\coordinate[label=below:C](C)at(1,0);
			\coordinate[label=above:D](D)at({(\d-2)/\d},1);
			\fill[lightgray](A)--(B)--(C)--(D)--cycle;
			\draw[dashed](1,0)node[above right]{$1$}--(1,1)--(0,1)node[left]{$1$};
			\draw[thick](A)--(B)--(C)--(D)--cycle;
			\draw($(A)!.5!(B)+(-.01,-.01)$)to[out=-140,in=0]+(-.2,-.1)node[left]{$\dfrac{2}{p}+\dfrac{d}{q}=d-2\sigma$};
			\draw($(C)!.5!(D)+(.01,.01)$)to[out=40,in=-90]+(.1,.2)node[above]{$\dfrac{2}{p}+\dfrac{d}{q}=d$};
			\foreach\P in{A,B,C,D}\fill(\P)circle[radius=\r];
		\end{tikzpicture}\caption{}\end{center}\end{figure}

For the singular value randomization, we have:
\begin{th.}\label{prop:a.s.linear}
	Let $d \ge 1$ and $0 \le \si < \frac{d}{2}$.
	Let $p,q \in [1,\I)$ satisfy $\frac{2}{p}+\frac{d}{q} = d-\si$ and $(\frac{1}{q}, \frac{1}{p}) \in ABCD$.
	When $d=2$, we omit the segment $AD$ (including the endpoints $A,D$) from $ABCD$.
	Let $\al = \min(p,q,2)$ and $r \in[\max(p,q),\infty)$.
	Then for any $\ga_0 \in \CH^{\si,\al}$, it holds that
	\begin{align}
		\|\sd^\si \rh(U(t)_\st \ga_0^{\om;\si})\|_{L^r_\om(\Om; L^p_t(\R;L^q_x))}
		\ls r^\tw \|\ga_0\|_{\CH^{\si,\al}}.
	\end{align}
\end{th.}

By using full randomization, we can control bigger exponents:
\begin{th.}\label{prop:a.s.linear2}
	Let $d\ge 1$ and $0 \le \si < \frac{d}{2}$.
	Let $p \in[2,\I)$ and $q\in [1,\infty)$ satisfy $\frac{2}{p}+\frac{d}{q}=d$, $\wh{q}\in[\max(q,2),\infty)$ and $r\in[\max(p,\wh{q}),\infty)$.
	Then for any $\gamma_0\in\CH^\si$, it holds that
	\begin{align}
		\|\sd^\si \rho(U(t)_\st \ga_0^{\om,\wo;\si})\|_{L^r_{\om,\wo}(\Om\times\wt{\Om}; L_t^p(\R;L_x^{\wh{q}}))}
		\lesssim r^{\frac{3}{2}}\|\ga_0\|_{\CH^\si}.
	\end{align}
\end{th.}

\begin{rmk.}
	Let $1 < q < \frac{d+1}{d-1}$ and $\frac{2}{p} + \frac{d}{q} = d$.
	On the one hand, by Theorem \ref{prop:a.s.linear}, we have
	\begin{align}
	\|\rh(U(t)_\st \ga_0^\om)\|_{L^p_t(\BR;L^q_x)} < \I \text{ for almost all } \om \in \Om
	\end{align}
	for any $\ga_0 \in \FS^{\min(p,q,2)}$.
	On the other hand, by Theorem \ref{SS}, we have
	\begin{align}
		\|\rh(U(t)_\st \ga_0^\om)\|_{L^p_t(\BR;L^q_x)} \ls \|\ga_0^\om\|_{\FS^\be},
	\end{align}
	where $\be := 2q/(q+1)$.
	It is easy to check that $\min(p,q,2) > \be$ always holds.
	In this sense, Theorem \ref{prop:a.s.linear} enables us to define the density function $\rh(U(t)_\st \ga_0)$ with $\ga_0 \in \FS^\al$, where $\al$ is bigger than the sharp exponent of the deterministic Strichartz estimates.
	We can say the same thing about Theorem \ref{prop:a.s.linear2}.
\end{rmk.}

\begin{rmk.}
The optimality of the Schatten exponents $\al$ in Theorems \ref{prop:a.s.linear} and \ref{prop:a.s.linear2} is a future problem.
\end{rmk.}

\begin{rmk.}
	The proofs of Theorems \ref{prop:a.s.linear} and \ref{prop:a.s.linear2} are quite easy and simple.
	However, to our best knowledge, these theorems are the first results dealing with this kind of randomization.
	We should emphasize that the orthogonality of eigenfunctions of $\ga_0$ is not essential in the proofs.
	In fact, we can prove
	\begin{align}
		&\No{\sum_{n=1}^\I a_n g_n(\om)|U(t)f_n|^2}_{L^r_\om L^p_t L^q_x} \ls r^\tw \|a\|_{\ell^\al}, \\
		&\No{\sum_{n=1}^\I a_n g_n(\om)|U(t)f_n^\wo|^2 }_{L^r_{\om,\wo} L^p_t L^q_x} \ls r^{\frac{3}{2}}\|a\|_{\ell^2},
	\end{align}
	under the same assumptions for exponents as in Theorems \ref{prop:a.s.linear} and \ref{prop:a.s.linear2} \textit{without orthogonality of }$(f_n)_{n=1}^\I$.
	This is because randomization works like orthogonality. (See Lemma \ref{lem:LDE}.)
\end{rmk.}

\subsubsection{Applications to the Hartree equation for infinitely many particles}
First, we give some applications to the local well-posedness of the Cauchy problem \eqref{eq:fNLH}.
The following results are not just applications of Theorems \ref{prop:a.s.linear} and \ref{prop:a.s.linear2};
we need some new ideas and tools to prove them
(see Lemmas \ref{lem:1D}, \ref{lem:2D}, \ref{lem:3D} and remarks to them).

\begin{th.}[One-dimensional case]\label{prop:Q 1,2}
	Let $d=1$.
	Assume that $w$ is a finite signed Borel measure on $\BR$, $f \in L^1_\xi \cap L^\I_\xi$, and $Q_0 \in \FS^2$ is self-adjoint.
	Then for almost all $\om \in \Om$, there exist $T>0$
	and a unique solution $Q(t) \in C(0,T;\FS^2)$ to
\begin{equation} \label{eq:random Q}
	\left\{ \,
	\begin{aligned}
		&i\pl_t Q = [-\De+w \ast \rh_Q,Q+\ga_f] \\
		&Q(0)=Q_0^\om
	\end{aligned}
	\right.
\end{equation} 
	such that $\rh_Q \in L^4_t(0,T;L^2_x)$.
\end{th.}

\begin{th.}[Two-dimensional case] \label{prop:Q 2}
	Let $d = 2$.
	Let $w \in L^1_x \cap L^{1+}_x$, $f \in L^1_\xi \cap L^\I_\xi$, and $Q_0$ be self-adjoint.
	Assume that $Q_0 \in \FS^{2-}$ or $Q_0 \in \CH^{0+}$ hold.
	Then for almost all $\om \in \Om$,
	there exist $T>0$
	and a unique solution $Q(t) \in C(0,T;\FS^2)$ to \eqref{eq:random Q} such that $\rh_Q \in L^2_t(0,T;L^2_x)$.
\end{th.}

\begin{rmk.}
	In the above statement, $w \in L^{1+}_x$ means there exists $\ep>0$ such that $w \in L^{1+\ep}_x$.
	Similarly, $Q_0 \in \FS^{2-}$ means there exists $\ep>0$ such that $Q_0 \in \FS^{2-\ep}$,
	and $Q_0 \in \CH^{0+}$ means there exists $\ep>0$ such that $Q_0 \in \CH^{\ep}$.
\end{rmk.}

\begin{th.}[Three-dimensional case] \label{prop:Q 3}
	Let $d = 3$.
	Let $w \in L^1_x \cap L^\I_x$, $f \in L^1_\xi \cap L^\I_\xi$, and $Q_0$ is self-adjoint.
	Assume that $Q_0 \in \FS^{3/2-}$ or $Q_0 \in \CH^{0+,3/2}$.
	Then for almost all $\om \in \Om$,
	there exist $T>0$
	and a unique solution $Q(t) \in C(0,T;\FS^2)$ to \eqref{eq:random Q} such that $w \ast \rh_Q \in L^2_t(0,T;L^2_x \cap L^\I_x)$.
\end{th.}

\begin{rmk.}
	Our results only claim the unique existence of local solutions; however, 
	we cannot construct any solution in our setting if we only use Theorem \ref{SS}.
We emphasize that our result allows strictly bigger Schatten exponents than Theorem \ref{SS}.
	\begin{table}[h]
		\begin{tabular}{|c|l|l|l|}
			\hline
			& $d=1$ & $d=2$ & $d=3$ \\ \hline
		\begin{tabular}{c} The Schatten exponents \\ that allowed by \\deterministic Strichartz estimates\end{tabular} & $\al < 2$ & $\al < \frac{3}{2}$ & $\al < \frac{4}{3}$ \\ \hline
		\begin{tabular}{c}The Schatten exponents\\that our results allow \end{tabular}
		&$\al =2$& $\al<2$& $\al < \frac{3}{2}$\\ \hline
		\end{tabular}
	\end{table}
\end{rmk.}

\begin{rmk.}
	To show Theorems \ref{prop:Q 1,2}, \ref{prop:Q 2}, and \ref{prop:Q 3}, we prove Lemmas \ref{lem:1D}, \ref{lem:2D}, and \ref{lem:3D}.
	These lemmas are similar to \cite[Theorems 5, 6]{2015LS}, but they have some advantages and we need new ideas and tools to prove them (see Remark \ref{rmk:compare}) 
\end{rmk.}
In Theorem \ref{prop:Q 2}, we need $\ep>0$ regularity for initial data $Q_0 \in \FS^2$.
However, by using full randomization, we can deal with all Hilbert--Schmidt initial data without any extra regularities:
\begin{cor.} \label{cor:full random 2D}
Let $d=2$, $f \in L^1_\xi \cap L^\I_\xi$, $w \in L^1_x \cap L^\I_x$, and $Q_0 \in \FS^2$ be self-adjoint.
Then for almost all $(\om,\wo) \in \Om\times\wt{\Om}$, there exists $T>0$ and a unique solution $Q(t) \in C(0,T;\FS^2)$ to
\begin{equation} \label{eq:full random Q}
	\left\{ \,
	\begin{aligned}
		&i\pl_t Q = [-\De+w \ast \rh_Q,Q+\ga_f] \\
		&Q(0)=Q_0^{\om,\wo}
	\end{aligned}
	\right.
\end{equation} 
such that $w \ast \rh_Q \in L^4_t(0,T;L^2_x \cap L^\I_x)$.	
\end{cor.}

\begin{rmk.}
	Corollary \ref{cor:full random 2D} immediately follows from Theorem \ref{prop:a.s.linear2} and \cite[Theorems 5]{2015LS}.
\end{rmk.}

Next, we give some application to the scattering problem of the linearized equation of \eqref{eq:fNLH}:
\begin{equation} \label{eq:L fNLH}
	\left\{ \,
	\begin{aligned}
		&i\pl_t Q = [-\De,Q]+[w \ast \rh_Q, \ga_f], \quad Q:\BR \to \CB(L^2_x), \\
		&Q(0)=Q_0.
	\end{aligned}
	\right.
\end{equation}
For given $f$ and $w$, we define a linear operator $\CL_1=\CL_1(f,w)$ by
\begin{align}
	\CL_1[g](t) := \rh \K{i\int_0^t U(t-\ta)_\st [w \ast g(\ta), \ga_f] d\ta}.
\end{align}

\begin{th.} \label{cor:singular scattering}
	Let $d \ge 2$ and $\si = \frac{d}{2}-1$.
	Let $\lxr^{2\si} f(\xi) \in L^1_\xi \cap L^\I_\xi$ and $w$ be a finite signed Borel measure on $\R^d$.
	Assume that $Q_0 \in \CH^{\si}$ is self-adjoint, and that $(1+\CL_1(f,w))^{-1} \in \CB(L^2_t(\BR_+;L^2_x))$. 
	Then for almost all $\om \in \Om $,
	there exists a unique global solution $Q(t) \in C(\BR_+;\CH^\si)$ to
	\begin{equation} \label{eq:singful random LQ}
		\left\{ \,
		\begin{aligned}
			&i\pl_t Q = [-\De,Q]+[w \ast \rh_Q, \ga_f], \\
			&Q(0)=Q_0^{\om;\si}
		\end{aligned}
		\right.
	\end{equation}
	such that $\rh_Q \in L^2_t(\BR_+;H^\si_x)$.
	Moreover, $Q(t)$ scatters; that is, there exists $Q_+ \in \FS^{\frac{2d}{d-1}}$ such that
	\begin{align}
		U(-t)Q(t)U(t) \to Q_+ \mbox{ in } \FS^{\frac{2d}{d-1}} \mbox{ as } t \to \I.
	\end{align}
\end{th.}

\begin{rmk.}
	Some sufficient conditions for $(1+\CL_1)^{-1} \in \CB(L^2_{t,x})$ are known.
	For example, see \cite[Corollary 1]{2014LS}, \cite[Propositions 1.5, 1.6]{2023Ha}.
\end{rmk.}

In Theorem \ref{cor:singular scattering}, we need $\si=\frac{d}{2}-1$ regularity for initial data.
However, we can omit this extra regularity by using full randomization:
\begin{th.} \label{cor:scattering}
	Let $d \ge 2$ and $\si = \frac{d}{2}-1$.
	Let $\lxr^{\si} f(\xi) \in L^1_\xi \cap L^\I_\xi$ and $w$ be a finite signed Borel measure on $\R^d$.
	Assume that $Q_0 \in \FS^2$ and $(1+\CL_1(f,w))^{-1} \in \CB(L^2_t(\BR_+;L^2_x))$. 
	Then for almost all $(\om,\wo) \in \Om\times\wt{\Om}$,
	there exists a unique global solution $Q(t) \in C(\BR_+;\FS^2)$ to
	\begin{equation} \label{eq:full random LQ}
		\left\{ \,
		\begin{aligned}
			&i\pl_t Q = [-\De,Q]+[w \ast \rh_Q, \ga_f], \\
			&Q(0)=Q_0^{\om,\wo}
		\end{aligned}
		\right.
	\end{equation}
	such that $\rh_Q \in L^2_t(\BR;L^2_x)$.
	Moreover, $Q(t)$ scatters; that is, there exists $Q_\pm \in \FS^{\frac{2d}{d-1}}$ such that
	\begin{align}
		U(-t)Q(t)U(t) \to Q_\pm \mbox{ in } \FS^{\frac{2d}{d-1}} \mbox{ as } t \to \I.
	\end{align}
\end{th.}

\begin{rmk.}
	Applications of singular value or full randomization to the nonlinear scattering problems of \eqref{eq:fNLH} are  interesting future problems. 
\end{rmk.}

\subsection{Strategy of the proof}
In this section, we explain the idea of the proof of Theorem \ref{prop:Q 1,2}, \ref{prop:Q 2}, and \ref{prop:Q 3}.
By Theorems \ref{prop:a.s.linear} and \ref{prop:a.s.linear2}, it suffices to prove the following lemmas:
\begin{lem.}[One-dimensional case] \label{lem:1D}
	Let $d=1$.
	Let $w$ be a finite signed Borel measure on $\BR$, $f \in L^1_\xi \cap L^\I_\xi$, and $Q_0$ be self-adjoint.
	Assume that $Q_0 \in \FS^2$ and $\|\rh(U(t)_\st Q_0)\|_{L^4_{t,\loc}(\R;L^2_x)} < \I$.
	Then there exist $T>0$ and a unique solution $Q(t) \in C(0,T;\FS^2)$ to \eqref{eq:fNLH}
	such that $\rh_Q \in L^4_t(0,T;L^2_x)$.
\end{lem.}

\begin{lem.}[Two-dimensional case] \label{lem:2D}
	Let $d=2$.
	Let $f \in L^1_\xi \cap L^\I_\xi$ and $Q_0$ be self-adjoint.
	Assume that there exists $\ep >0$ such that $w \in L^1_x \cap L^{1+\ep}_x$.
	Assume that $Q_0 \in \FS^{2-\ep}$ or $Q_0 \in \CH^{\ep}$ hold, and $\|\rh(U(t)_\st Q_0)\|_{L^2_{t,\loc}(\BR;L^2_x)} < \I$.
	Then there exist $T>0$ and a unique solution $Q(t) \in C(0,T;\FS^2)$ to \eqref{eq:fNLH}
	such that $\rh_Q \in L^2_t(0,T;L^2_x)$.
\end{lem.}

\begin{lem.}[Three-dimensional case]\label{lem:3D}
	Let $d = 3$.
	Let $w \in L^1_x \cap L^\I_x$, $f(\xi) \in L^1_\xi \cap L^\I_\xi$, and $Q_0$ be self-adjoint.
	Assume that there exists $\ep>0$ such that $Q_0 \in \FS^{3/2-\ep}$ or $Q_0 \in \CH^{\ep,3/2}$.
	Assume that $\|w \ast \rh(U(t)_\st Q_0)\|_{L^2_{t,\loc}(\R;L^2_x \cap L^\I_x)} < \I$.
	Then there exist $T>0$ and a unique solution $Q(t) \in C(0,T;\FS^2)$ to \eqref{eq:fNLH}
	such that $w \ast \rh_Q \in L^2_t(0,T;L^2_x \cap L^\I_x)$.
\end{lem.}

\begin{rmk.} \label{rmk:compare}
	The above lemmas are similar to \cite[Theorems 5, 6]{2015LS}.
	One of these results never includes the other; hence, there is no simple superiority or inferiority.
	However, our results have some advantages:
	\begin{itemize}
		\item In \cite[Theorem 5]{2015LS}, we need $w \in L^1_x \cap L^\I_x$ to get local solution; however, all general finite signed measures are acceptable in Lemma \ref{lem:1D}, and $w \in L^1_x \cap L^{1+\ep}_x$ is sufficient in Lemma \ref{lem:2D}.
		\item In \cite[Theorem 5]{2015LS}, we only get information of the indirect object $w \ast \rh_Q$; however, we can have more direct information of $\rh_Q$ in Lemmas \ref{lem:1D} and \ref{lem:2D}.   
		\item In \cite[Theorem 6]{2015LS} with $d=3$, we need $\lxr f(\xi) \in L^2_x$ and $w \in W^{1,1}_x \cap W^{1,\I}_x$ and $\sd Q_0 \in \FS^2$. However, we do not need this type of regularities for $f$, $w$, and $Q_0$ in Lemma \ref{lem:3D}. 
	\end{itemize}
\end{rmk.}

To prove the above lemmas, we need some new ideas and tools.
We use the orthonormal Strichartz estimates improved by Bez, Hong, Lee, Nakamura, and Sawano in \cite{2019BHLNS}.
See Theorem \ref{th:Ostri}.
As far as the authors know, this is the first local-in-time result applying it.

In this paper, we prove the unique existence of the local solution by the standard contraction mapping argument.
We consider the following Duhamel integral for the system $(Q,\rh_Q)$:
\begin{align} \label{eq:D}
	\begin{dcases}
		&Q(t) = U(t)_\st Q_0 - i\int_0^t U(t-\ta)_\st [w \ast \rh_Q(\ta),Q(\ta)+\ga_f] d\ta, \\
		&\rh_Q(t)= \rh(U(t)_\st Q_0) - i\int_0^t \rh\Big(U(t-\ta)_\st [w \ast \rh_Q(\ta),Q(\ta)+\ga_f]\Big) d\ta.
	\end{dcases}
\end{align}
If we use the Minkowski inequality directly to estimate the Duhamel term, we get
\begin{align}
	\No{\int_0^T U(\ta)^*V(\ta)Q(\ta) U(\ta) d\ta}_{\FS^\al}
	\le \|V\|_{L^1_T L^\I_x} \|Q\|_{C_T \FS^\al}
\end{align}
and we cannot avoid the appearance of $\|V\|_{L^1_T L^\I_x}$.
To overcome this problem, in \cite{2015LS}, Lewin and Sabin used the density functional method;
that is, they reduced the original problem to the nonlinear equation of density function $\rh_Q$.
Later, they recovered the original solution $Q(t)$ from its density function $\rh_Q$.
However, this method is a little bit complicated because we need some nontrivial calculations in order to estimate two different solutions when we use the contraction mapping argument.
In this paper, we prove the key estimate Lemma \ref{lem:bilinear} instead of using density functional method.

To estimate a priori bound of $\rh_Q$, the above Duhamel integral \eqref{eq:D} is not enough to get sufficient Schatten exponents.
Therefore, we repeat the Duhamel integral and modify the set-up for contraction mapping argument. See \eqref{eq:Duhamel} and \eqref{eq:Duhame2}.

\section{Probabilistic Strichartz estimates}
In this section, we prove the main estimates in this paper.
\subsection{Basic tools}
In this section, we collect some necessary tools.
\begin{lem.}[\cite{2014LM}; Unit-scale Bernstein inequality]\label{lem:unit-scale B}
	For all $1\le r_1\le r_2\le\infty$ and for all $k\in\Z^d$, it holds that
	\EQ{
		\|\Pi_kf\|_{L^{r_2}(\R^d)}\le C(r_1,r_2)\|\Pi_kf\|_{L^{r_1}(\R^d)},
	}
	where $C(r_1,r_2)$ is a constant independent of $k\in\Z^d$.
\end{lem.}

\begin{lem.}[\cite{BT08a}; Large deviation estimate]\label{lem:LDE}
Let $(g_n)_{n=1}^\I$ be a sequence of random variables introduced in Section \ref{sec:randomization}.
Then there is a constant $C>0$ such that
\begin{align}
  \nor{\sum_{n=1}^\I a_n g_n}_{L^{r}_{\omega}(\Om)} \le C r^\tw \|a_n\|_{\ell^2_n}.
\end{align}
holds for any $r\in[2,\infty)$ and $(a_n)_n\in\ell^2_n$.
\end{lem.}

\begin{lem.}\label{lem:a.s.linear function}
	Let $d \ge 1$.
	Let $(p,q)\in[2,\infty)^2$ be an admissible pair of the standard Strichartz estimates.
	Let $\wh{q}\in[q,\infty)$ and $r\in[\max(p,\wh{q}),\infty)$.
	For any $f\in L^2(\R^d)$, we have
	\begin{align}
		&\|U(t)f^{\wo}\|_{L^r_{\wo}(\wt{\Om}; L_t^pL_x^\wq(\R\times\R^d))}
		\lesssim r^\tw\|f\|_{L_x^2(\R^d)}.
	\end{align}
\end{lem.}

\begin{proof}
	By the Minkowski inequality and Lemma \ref{lem:LDE}
	\begin{align}
		\|U(t)f^\wo\|_{L^r_{\wo}(\wt{\Omega}; L_t^pL_x^{\wh{q}})}
		&\le\Biggl\|\biggl\|\sum_{k\in\Z^d}U(t)h_k(\wo)\Pi_kf\biggr\|_{L^r_{\wo}(\wt{\Om})}\Biggr\|_{L_t^pL_x^{\wh{q}}} \\
		&\ls r^\tw \No{\|U(t)\Pi_kf\|_{\ell^2_k}}_{L_t^pL_x^{\wh{q}}}.
	\end{align}
	Combining the Minkowski inequality, Lemma \ref{lem:unit-scale B} and Strichartz estimate, we have
	\begin{align}
		&\|U(t)f^\wo\|_{L^r_{\wo}(\wt{\Omega}; L_t^pL_x^{\hat{q}})}
		\lesssim r^\tw \No{\|U(t)\Pi_kf\|_{L_t^pL_x^{\wh{q}}}}_{\ell^2_k}
		\\&\lesssim r^\tw \No{\|U(t)\Pi_kf\|_{L_t^pL_x^q}}_{\ell^2_k}
		\lesssim r^\tw \K{\sum_{k \in \Z^d}\|\Pi_kf\|_{L_x^2}^2}^{1/2}
		= r^\tw \|f\|_{L_x^2},
	\end{align}
	which concludes the desired estimate.
\end{proof}

\subsection{Proof of Theorem \ref{prop:a.s.linear}}
We write the singular value decomposition of $\sd^\si_\st \ga_0$
as $\sd^\si_\st \ga_0 =\sum_{n}a_n|u_n\rangle \langle v_n|$.
Then we have
\begin{align}
	\ga_0^{\om;\si} &= \sum_{n=1}^\I a_n g_n(\om)|\sd^{-\si} u_n\rangle \langle \sd^{-\si}v_n|
	=: \sum_{n=1}^\I a_n g_n(\om)|\wt{u_n} \rangle \langle \wt{v_n}|
\end{align}
By the Minkowski inequality and Lemma \ref{lem:LDE}, we have
	\begin{align*}\label{ineq:a.s.linear}
		&\|\sd^\si \rh(U(t)_\st\ga_0^{\om;\si})\|_{L^r_\om(\Om; L_t^pL_x^q)} \\
		&\quad =\nor{\sum_{n=1}^\I a_n g_n(\om) \sd^\si \K{\ov{U(t)\wt{v_n}} U(t)\wt{u_n}} }_{L^r_\omega(\Omega; L_t^pL_x^q)} \\
		&\quad \le\nor{\No{\sum_{n=1}^\I a_ng_n(\omega) \sd^\si \K{\ov{U(t)\wt{v_n}} U(t)\wt{u_n}} }_{L^r_\omega(\Omega)}}_{L_t^pL_x^q} \\
		&\quad \ls r^\tw\nor{\nor{a_n \sd^\si \K{\ov{U(t)\wt{v_n}} U(t)\wt{u_n}}}_{\ell_n^2}}_{L_t^pL_x^q}.
	\end{align*}
Since $(\frac{1}{q}, \frac{1}{p}) \in ABCD$, there exists $p_0,q_0,p_1,q_1 \in [1,\I]$ such that
\begin{align}
	\frac{2}{p_0}+\frac{d}{q_0}=d, \quad \frac{2}{p_1}+\frac{d}{q_1}=d-2\si, \quad \tw\K{\frac{1}{p_0} + \frac{1}{p_1}} = \frac{1}{p}, \quad \tw\K{\frac{1}{q_0} + \frac{1}{q_1}} = \frac{1}{q}.
\end{align}
Using $\ell^{\alpha}\subset\ell^2$, the Minkowski inequality, the Strichartz estimate, and the fractional Leibniz rule, we obtain
	\begin{align}
		\nor{\nor{a_n \sd^\si \K{\ov{U(t)\wt{v_n}} U(t)\wt{u_n}}}_{\ell_n^2}}_{L_t^pL_x^q}
		&\le\nor{\nor{a_n \sd^\si \K{\ov{U(t)\wt{v_n}} U(t)\wt{u_n}}}_{\ell_n^{\alpha}}}_{L_t^pL_x^q}\\
		&\le\nor{\nor{a_n \sd^\si \K{\ov{U(t)\wt{v_n}} U(t)\wt{u_n}}}_{L_t^pL_x^q}}_{\ell_n^{\alpha}} \\
		&\ls\No{a_n
		\|U(t)\wt{u_n}\|_{L_t^{2p_1}L_x^{2q_1}}
		\|U(t)\wt{v_n}\|_{L_t^{2p_0}H_x^{\si,2q_0}} }_{\ell^\al_n} \\
		&\quad +
		\No{a_n\|U(t)\wt{u_n}\|_{L_t^{2p_0}H_x^{\si,2q_0}} 
		\|U(t)\wt{v_n}\|_{L_t^{2p_1} L_x^{2q_1}}}_{\ell^\al_n} \\
		&\ls \|a\|_{\ell^\al} = \|\ga_0\|_{\CH^{\si,\al}}.
 	\end{align}
	\qed

\subsection{Proof of Theorem \ref{prop:a.s.linear2}}
We write the singular value decomposition of $\sd^\si_\st \ga_0$ as $\sd^\si_\st\ga_0=\sum_{n}a_n|u_n\rangle \langle v_n|$.
Then we have
\begin{align}
	\ga_0^{\om,\wo;\si} &= \sum_{n=1}^\I a_n g_n(\om)|\sd^{-\si} u_n^\wo \rangle \langle \sd^{-\si}v_n^\wo|
	=: \sum_{n=1}^\I a_n g_n(\om)|\wt{u_n}^\wo \rangle \langle \wt{v_n}^\wo|
\end{align}
By the Minkowski inequality, Lemmas \ref{lem:LDE} and \ref{lem:a.s.linear function}, we have
	\begin{align}
	&\|\sd^\si \rho(U(t)_\st \ga_0^{\om,\wt{\om}})\|_{L^r_{\om,\wo}(\Om\times\wt{\Om}; L_t^p L_x^{\wh{q}})}  \\
	&\quad = \No{\sum_{n=1}^\I a_n g_n(\om) \sd^\si\K{\ov{U(t)\wt{v_n}^\wo}U(t)\wt{u_n}^\wo } }_{L^r_{\om,\wo} L^p_t L^\wq_x} \\
	&\quad \ls r^\tw \No{ a_n \sd^\si\K{\ov{U(t)\wt{v_n}^\wo}U(t)\wt{u_n}^\wo } }_{L^r_\wo L^p_t L^\wq_x \ell_n^2} \\
	&\quad \le r^\tw \No{ a_n \sd^\si\K{\ov{U(t)\wt{v_n}^\wo}U(t)\wt{u_n}^\wo } }_{\ell_n^2 L^r_\wo L^p_t L^\wq_x } \\
	&\quad \ls r^\tw \No{ a_n \|U(t)u_n^\wo\|_{L^{2r}_\wo L^{2p}_t L^{2\wq}_x}  \|U(t)v_n^\wo\|_{L^{2r}_\wo L^{2p}_t L^{2\wq}_x}  }_{\ell_n^2} \ls r^{\frac{3}{2}}\|\gamma_0\|_{\mathfrak{S}^2}.
\end{align}
\qed

\subsection{No gain of regularity by the randomizations}

We can obtain neither lower Schatten integrability nor higher Sobolev regularity by the singular value randomization unless $(g_n)_{n=1}^\I$ doesn't gather to 0.
This argument is similar to \cite[Lemma B.1]{BT08a}.

\begin{prop.}\label{prop:appendix lower}
Assume that the sequence of random variables $(g_n)_{n=1}^\I$ introduced in Section \ref{sec:randomization} satisfies the following condition: there exists $c>0$ such that
\EQ{\label{ineq:appendix lower}
  \limsup_{n\to\infty}\mathbb{P}(\{\omega\in\Omega : |g_n(\omega)|\le c \})<1.
}
Let $1 \le \be \le \al < \I$ and $\si \in \BR$.
Then $\sd^\si_\st \ga_0 \in \FS^{\al} \setminus \FS^{\be}$ implies $\sd^\si_\st \ga_0^{\om;\si} \in \FS^\al \setminus \FS^\be$ for almost all $\om \in \Om$.
\end{prop.}

\begin{rmk.}
If $(g_n)_{n=1}^\I$ is a sequence of Gaussian random variables with uniformly lower bounded variances, \eqref{ineq:appendix lower} is satisfied.
\end{rmk.}

\begin{proof}
We can assume that $\si = 0$.
First, we prove $\ga_0^\om \in \FS^\al$ for almost all $\om \in \Om$.
We have by \eqref{eq:distr g_n}
\begin{align}
	\BE \Dk{\sum_{n=1}^\I |a_n g_n(\om)|^\al} = \sum_{n=1}^{\I} |a_n|^\al \BE[|g_n|^\al] \ls \|a\|_{\ell^\al}^{\al} < \I.
\end{align}
Therefore, we obtain
\begin{align}
	\sum_{n=1}^\I |a_n g_n(\om)|^\al < \I
\end{align}
for almost all $\om \in \Om$.
	
Next, we prove $\ga_0^\om \notin \FS^\be$ for almost all $\om \in \Om$.
From the assumption \eqref{ineq:appendix lower}, there exists $\delta\in(0,1)$ and $N \in \BN$ such that $\mu_n([-c,c])\le1-\delta$ holds for any $n \ge N$. Hence,
\begin{align}
  &\int_{\Omega}e^{-\|\gamma_0^\omega\|_{\mathfrak{S}^{\beta}}^{\beta}}\,d\mathbb{P}(\omega)
  =\int_{\Omega}e^{-\sum_{n=1}^{\infty}|a_ng_n(\omega)|^{\beta}}\,d\mathbb{P}(\omega)
  \\&=\prod_{n=1}^{\infty}\K{\int_{\{|x|\le c\}}e^{-|a_nx|^{\beta}}\,d\mu_n(x)+\int_{\{|x|> c\}}e^{-|a_nx|^{\beta}}\,d\mu_n(x)}
  \\&\le\prod_{n=N}^{\infty}\K{\mu_n([-c,c])(1-e^{-|a_nc|^{\beta}})+e^{-|a_nc|^{\beta}}}
  \\&\le\prod_{n=N}^{\infty}\K{(1-\delta)+e^{-|a_nc|^{\beta}}\delta}
  =\prod_{n=N}^{\infty}\K{1-\delta(1-e^{-|a_nc|^{\beta}})}.
\end{align}
Since $|a_n c|^{\beta}\to+0$ from $\gamma_0\in\mathfrak{S}^{\alpha}$, there exists $\lambda\in(0,1)$ such that
\EQ{
  1-e^{-|a_nc|^{\beta}}\ge\lambda|a_nc|^{\beta}
}
for all $n\in\N$.
Combining $\|\gamma_0\|_{\mathfrak{S}^{\beta}}=\infty$,
\EQ{
  \sum_{n=N}^{\infty}(1-e^{-|a_nc|^{\beta}})
  \ge\lambda c^{\beta}\sum_{n=N}^{\infty}|a_n|^{\beta}=\infty.
}
Therefore, we obtain that
\begin{align}
  \int_{\Omega}e^{-\|\gamma_0^\omega\|_{\mathfrak{S}^{\beta}}^{\beta}}\,d\mathbb{P}(\omega)
  &\le\prod_{n=N}^{\infty}\K{1-\delta(1-e^{-|a_nc|^{\beta}})}
  \\&\le\exp{\biggl(-\delta\sum_{n=N}^{\infty}(1-e^{-|a_nc|^{\beta}})\biggr)}
  =0,
\end{align}
which concludes that $\|\gamma_0^\omega\|_{\mathfrak{S}^{\beta}}=\infty$ almost surely.
\end{proof}

We also cannot gain any regularities by the full randomization.

\begin{prop.} \label{prop:no gain}
	Let $(g_n)_{n=1}^\I$ introduced in Section \ref{sec:randomization} satisfies \eqref{ineq:appendix lower} for $c>0$.
	Moreover, assume that $(h_k)_{k\in\Z^d}$ is i.i.d. of the Rademacher distribution. 
	Let $1 \le \be \le \al < \I$ and $\si \in \BR$.
	Then $\sd^\si_\st \ga_0 \in \FS^{\al} \setminus \FS^{\be}$ implies $\sd^\si_\st \ga_0^{\om,\wo;\si} \in \FS^\al \setminus \FS^\be$ for almost all $(\om,\wo) \in \Om\times\wt{\Om}$.
\end{prop.}

\begin{proof}
We can assume that $\si = 0$.
By Proposition \ref{prop:appendix lower}, it suffices to show that $\CR^\wo \ga_0 \CR^\wo$ belongs to $\FS^\al \setminus \FS^\be$ for almost all $\wo\in \wt{\Om}$.
We can write $\CR^\wo = \CF^{-1} R^\wo(\xi) \CF$, where 
\begin{equation}
	R^\wo(\xi)= \sum_{k\in\Z^d} h_k(\wo)\ch_k(\xi).
\end{equation}
Since the overlap of $(\supp \ch_k)_{k\in \Z^d}$ is a null set and $|h_k(\wo)|\equiv1$ for almost all $\wo \in \wt{\Om}$, we have $|R^\wo(\xi)|=1$ for almost all $(\wo,\xi) \in \wt{\Om}\times \R^d$.
Therefore, we have
\begin{equation}
	\No{\CR^\wo \ga_0 \CR^\wo}_{\FS^p}= \No{\ga_0}_{\FS^p}
\end{equation}
for all $p \in [1,\I]$.
\end{proof}

\section{Application to the Hartree equation: Unique existence of the local solution}
In this section, we prove Lemmas \ref{lem:1D}, \ref{lem:2D}, and \ref{lem:3D}.
\subsection{Preliminaries}
One of the most important tools is the orthonormal Strichartz estimates:
\begin{th.}[\cite{2014FLLS}Theorem 1; \cite{2017FS}Theorem 8; \cite{2019BHLNS}Theorem 1.5] \label{th:Ostri}
	Let $d \ge 1$.
	Let $(p,q,\al)$ be admissible and $\frac{2}{p}+\frac{d}{q}=d-2s$ for $s \in (0,d/2)$,
	where we call $(p,q,\al)$ admissible when $\frac{1}{\al} \ge \frac{1}{dp} + \frac{1}{q}$ and $\al < p$.
	Then it holds that
	\begin{align}
		\|\rh(U(t)_\st Q_0)\|_{L^p_t L^q_x} \ls \|\sd^s_\st Q_0\|_{\FS^\al}.
	\end{align}
\end{th.}

Let $H$ and $K$ be separable complex Hilbert spaces.
For $\al \in [1,\I]$ and compact operator $A: H \to K$,
we define Schatten $\al$--norm by
\begin{align}
	\|A\|_{\FS^\al(H\to K)}
	:= 
	\begin{cases}
		&\K{\Tr_H(|A^*A|^{\frac{\al}{2}})}^{\frac{1}{\al}} \mbox{ if } 1 \le \al <\I,\\
		&\|A\|_{\CB(H \to K)} \mbox{ if } \al = \I,
	\end{cases}
\end{align}
where $\Tr_H$ is the trace in $H$ and $\CB(H\to K)$ is the space of all bound linear operators from $H$ to $K$.
For each $\al \in [1, \I]$, we write
\begin{align}
	&\FS^\alpha = \FS^\alpha(L^2_x \to L^2_x),
	\quad \FS^\al_{t,x} := \FS^\al(L^2_{t,x} \to L^2_{t,x}), \\
	&\FS^\al_{x \to\tx} := \FS^\al(L^2_x \to L^2_{t,x}),
	\quad \FS^\alpha_{\tx \to x} := \FS^\al(L^2_{t,x} \to L^2_x).
\end{align}
Now, we recall the notations, definitions, and results in \cite[Section 3]{2023H}.
Let $I \subset \R$ be an interval.
Assume that $A(t) \in \CB(L^2_x)$ for all $t \in I$, $\sup_{t \in I} \|A(t)\|_{\CB} < \I$,
and $I \ni t \mapsto A(t) \in \CB(L^2_x)$ is strongly continuous.
Define $A^\ocircle : L^2_x \to C_t(I;L^2_x)$ and $A^\oast : L^1_t(I;L^2_x) \to L^2_x$ by
\begin{align}
	&(A^\ocircle u_0)(t,x) := (A(t)u_0)(x), \\
	&(A^\oast g)(x) := \int_I A(\ta)^* g(\ta,x) d\ta.
\end{align}
Note that $A^\ocircle$ and $A^\oast$ are formally adjoint to each other.
The following lemma is useful.
\begin{lem.}[Duality principle; \cite{2017FS} Lemma 3; see also \cite{2023H} Lemma 3.1] \label{lem:DP}
	Let $p,q,\al \in [1,\I]$.
	The followings are equivalent:
	\begin{itemize}
	\item For any $\ga_0 \in \FS^\al$, 
	\begin{align}
		\|\rh(A(t)_\st \ga_0)\|_{L^p_t(I;L^q_x)} \le C \|\ga_0\|_{\FS^\al}.
	\end{align}
	
	\item For any $g \in L^{2p'}_t(I;L^{2q'}_x)$, 
	\begin{align}
		\|g A^\oc\|_{\FS^{2\al'}_{x \to \tx}} \le C'\|g\|_{L^{2p'}_t(I;L^{2q'}_x)}.
	\end{align}
	
	\item For any $g \in L^{2p'}_t(I;L^{2q'}_x)$, 
	\begin{align}
		\|A^\os g\|_{\FS^{2\al'}_{\tx \to x}} \le C''\|g\|_{L^{2p'}_t(I;L^{2q'}_x)}.
	\end{align}
	\end{itemize}
	Moreover, $\sqrt{C}$, $C'$ and $C''$ coincide.
\end{lem.}

The following Chirst--Kiselev type lemma in Schatten classes is useful.
Note that the following result is already obtained in \cite{1970Go-Kre} before \cite{2001CK}.
Therefore, we should call it Gohberg--Kre\u{\i}n theorem.

\begin{th.}[Gohberg--Kre\u{\i}n; see \cite{2003BS}, \cite{1970Go-Kre}; see also \cite{2023H}] \label{SCK}
	Let $d \ge 1$, $T \in (0,\I]$ and $\al \in (1,\I)$.
	Suppose that we can write $\BT \in \CB(L^2_{t,x})$ by
	\begin{align}
		(\BT g)(t,x) = \int_0^T K(t,\ta)g(\ta)d\ta.
	\end{align}
	If $\BT \in \FS^\al(L^2_{t,x})$, then  $\BD$ defined by
	\begin{align}
		(\BD g)(t,x) = \int_0^{t} K(t,\ta)g(\ta)d\ta
	\end{align}
	is also in $\FS^\al(L^2_{t,x})$, and we have
	\begin{align}
		\|\BD\|_{\FS^\al(L^2_{t,x})} \le C_\al \|\BT\|_{\FS^\al(L^2_{t,x})}.
	\end{align}
\end{th.}

We also use the standard Christ--Kiselev theorem.
We can prove the following theorem by the standard argument with the Whitney decomposition.
\begin{th.}[Christ--Kiselev; See \cite{2001CK}] \label{th:OCK}
Let $T \in (0,\I]$.
Let $X_1$, $X_2$ and $X$ be Banach spaces.
	Assume that bilinear maps $B,B_<:L^{p_1}_t(0,T;X_1) \times L^{p_2}_t(0,T;X_2) \to L^p_t(0,T;X)$ are defined by
	\begin{align}
		&B(g_1,g_2):= \int_0^T dt_1 K_1(t,t_1)g_1(t_1) \int_0^T dt_2 K_2(t_1,t_2)g_2(t_2), \\
		&B_<(g_1,g_2):= \int_0^{t} dt_1 K_1(t,t_1)g_1(t_1) \int_0^{t_1} dt_2 K_2(t_1,t_2)g_2(t_2).
	\end{align}
	If $p_1<p$ and $1/p_1+1/p_2 > 1$, then we have
	\begin{align}
		&\|B(g_1,g_2)\|_{L^p_t(0,T;X)}
		\ls \|g_1\|_{L^{p_1}_t(0,T;X_1)}
		\|g_2\|_{L^{p_2}_t(0,T;X_2)} \\
		&\quad \Longrightarrow
		\|B_<(g_1,g_2)\|_{L^p_t(0,T;X)}
		\ls \|g_1\|_{L^{p_1}_t(0,T;X_1)}
		\|g_2\|_{L^{p_2}_t(0,T;X_2)}.
	\end{align}
\end{th.}

\subsection{Key estimate}
\begin{lem.}[Key estimate] \label{lem:bilinear}
	Let $d \ge 1$ and $2 \le \al \le \I$.
	Let $T \in (0,\I]$.
	Assume that $\mu,\nu \in[1,\I]$ satisfy $\frac{2}{\mu}+\frac{d}{\nu}=2$ and
	\begin{itemize}
		\item $1 \le \mu \le \frac{4}{3}$ when $d=1$,
		\item $1 \le \mu < 2$ when $d=2$,
		\item $1 \le \mu \le 2$ when $d\ge 3$.
	\end{itemize}
	Then we have
	\begin{align*}
		\No{\int_0^T U(\ta)^* V(\ta) Q(\ta) U(\ta) d\ta }_{\FS^{\al}}
		\lesssim \|V\|_{L^\mu_t(0,T;L^\nu_x)}\|Q\|_{C(0,T;\FS^{\al})}.
	\end{align*}
\end{lem.}
\begin{proof}
	We assume $d \ge3$ since we can prove the other cases similarly.
	It is straightforward to verify the case $\mu=1$.
	Hence, we assume that $\mu =2$.
	First, we consider the case $\al = 2$.
	Let $\frac{1}{q}:= \tw - \frac{1}{d}$.
	Then we have by the endpoint Strichartz estimate (see \cite{1998KT})
	\begin{align*}
		\abs{ \Tr \Dk{ A \int_0^T U(\ta)^* V(\ta) Q(\ta) U(\ta) d\ta } } 
		&= \abs{ \int_0^T d\ta \Tr \Dk{ Q(\ta) U(\ta) A U(\ta)^* V(\ta)  } } \\
		&\le \No{\rh\Dk{Q(t)(U(t)_\st A)}}_{L^2_T L^{d'}_x}
		\|V\|_{L^2_T L^d_x} \\
		&\le \|\rh_{Q^2}^{1/2}\|_{C_T L^2_x}
		\|\rh(U(t)_\st |A|^2)^\tw\|_{L^2_T L^{q}_x}
		\|V\|_{L^2_T L^d_x} \\
		&\le \|Q\|_{C_T \FS^2}
		\|A\|_{\FS^2}
		\|V\|_{L^2_T L^d_x}. 
	\end{align*}
	Next, we consider the case $\al=\I$.
	Let $A \in \FS^1$.
	Then we have
	\begin{align*}
		&\abs{ \Tr \Dk{ A \int_0^T U(\ta)^* V(\ta) Q(\ta) U(\ta) d\ta } } \\
		&\quad \le \No{\rh\Dk{Q(t)(U(t)_\st A)}}_{L^2_T L^{d'}_x}
		\|V\|_{L^2_T L^d_x} \\
		&\quad \le \|\rh(Q(t) U(t) |A| U(t)^* Q(t)^*)^\tw\|_{L^\I_T L^2_x}
		\|\rh(U(t)_\st |A|)^\tw \|_{L^2_T L^q_x}
		\|V\|_{L^2_T L^d_x} \\
		&\quad \le  \|Q\|_{C_T \FS^\I}
		\|A\|_{\FS^1}
		\|V\|_{L^2_T L^d_x}.
	\end{align*}
\end{proof}

\subsection{Proof of Lemma \ref{lem:1D}}
We have the Duhamel formula
\begin{align} \label{eq:Duhamel}
	\left \{
	\begin{aligned}
		Q(t) &= \Ph_1[Q,\rh_Q], \\
		\rh_Q(t) &= \Ph_2[Q,\rh_Q],
	\end{aligned}
	\right.
\end{align}
where
\begin{align}
	\Ph_1[Q,\rh_Q]&:= U(t)_\st Q_0 + \CD_V[Q] +\CD_V[\ga_f] \\
	&:=U(t)_\st Q_0 - i \int_0^t U(t-\ta)_\st [V(\ta),Q(\ta)] d\ta- i \int_0^t U(t-\ta)_\st [V(\ta),\ga_f] d\ta, \\
	V&:= w\ast \rh_Q, 
\end{align}
and
\begin{align}
	\Ph_2[Q,\rh_Q](t)
	&:=\rh(U(t)_\st Q_0) + \rh(\CD_V[U(t)_\st Q_0]) \\
	&\quad + \rh(\CD_V^2[Q]) + \rh(\CD_V^2[\ga_f]) - \CL_1(\rh_Q).
\end{align}
We used the notation $\CL_1[\rh_Q]:= -\rh(\CD_{w \ast \rh_Q}[\ga_f])$.
Define
\begin{align}
	E(T,R):= \{(Q,g) \in C_T \FS^2 \times L^4_T L^2_x : \|Q\|_{C_T \FS^2} \le R, \|g\|_{L^4_T L^2_x}\le R\},
\end{align}
where $R:= 2(\|Q_0\|_{\FS^2} + \|\rh(U(t)_\st Q_0)\|_{L^4_{T} L^2_x})$ and $T>0$ will be chosen later.
We prove that $\Ph : E(T,R) \ni (Q,g) \to (\Ph_1[Q,g],\Ph_2[Q,g]) \in E(T,R)$ is a contraction mapping.

Let $(Q,g) \in E(T,R)$.
Then by Lemma \ref{lem:bilinear} and Kato--Seiler--Simon's inequality (see \cite{1975SS} and \cite[Theorem 4.1]{2005Simon}), we have
\begin{align}
	\|\Ph_1[Q,g](t)\|_{\FS^2}
	&\le \|Q_0\|_{\FS^2}
	+ 2\No{\int_0^t U(\ta)^*V(\ta)Q(\ta)U(\ta) d\ta }_{\FS^2} \\
	&\quad + 2 \No{\int_0^t U(\ta)^*V(\ta)\ga_f(\ta)U(\ta) d\ta }_{\FS^2}\\
	&\le \|Q_0\|_{\FS^2} + C \|V\|_{L^{4/3}_T L^2_x} \|Q\|_{C_T \FS^2}
	+ 2\int_0^T \|V(\ta) \ga_f\|_{\FS^2} d\ta \\
	&\le \frac{R}{2} + C T^\tw R^2 + C \|f\|_{L^2_\xi} T^{\frac{1}{4}} R
	\le R
\end{align}
if we choose $T>0$ sufficiently small.
For the density function, we have
\begin{align}
	\|\Ph_2[Q,g]\|_{L^4_T L^2_x}
	&\le \|\rh(U(t)_\st Q_0)\|_{L^4_T L^2_x}
	 + \|\rh(\CD_V[U(t)_\st Q_0])\|_{L^4_T L^2_x} \\
	 &\quad + \|\rh(\CD_V^2[Q])\|_{L^4_T L^2_x}
	 + \|\rh(\CD_V^2[\ga_f])\|_{L^4_T L^2_x} + \|\CL_1(\rh_Q)\|_{L^4_T L^2_x} \\
	&=:A_1 + A_2 + A_3 + A_4+A_5.
\end{align}

\subsubsection{Estimate for $A_1$}
By definition of $R$, we obtain
\begin{align}
	A_1 \le \frac{R}{2}.
\end{align}

\subsubsection{Estimate for $A_2$}
Let
\begin{align}
	\CV:= \int_0^T U(\ta)^* V(\ta) U(\ta) d\ta.
\end{align}
Since it follows from Theorem \ref{SS} and the duality argument that
\begin{align}
	\No{\rh(U(t)[\CV,Q_0]U(t)^*)}_{L^4_T L^{2}_x}
&\ls 
	\|\CV Q_0 \|_{\FS^{4/3}} \ls \|V\|_{L^{4/3}_T L^2_x} \|Q_0\|_{\FS^2},
\end{align}
by the Christ--Kiselev lemma (see \cite{2001CK}), we obtain
\begin{align}
	A_2 \le C\|V\|_{L^{4/3}_T L^2_x} \|Q_0\|_{\FS^2} \le C T^\tw R^2 \le \frac{R}{8}
\end{align}
if we choose sufficiently small $T > 0$.

\subsubsection{Estimate for $A_3$}
Note that by Lemma \ref{lem:bilinear}
\begin{align}
	&\No{\rh \Big(U(t) \int_0^T U(\ta)^* V(\ta)U(\ta) d\ta \int_0^T U(\ta_1)^* V(\ta_1)Q(\ta_1) U(\ta_1) d\ta_1 U(t)^*\Big)}_{L^4_T L^2_x} \\
	&\quad \ls \No{ \int_0^T U(\ta)^* V(\ta)U(\ta) d\ta \int_0^T U(\ta_1)^* V(\ta_1)Q(\ta_1) U(\ta_1) d\ta_1 }_{\FS^{4/3}} \\
	&\quad \le \No{ \int_0^T U(\ta)^* V(\ta)U(\ta) d\ta }_{\FS^4}
	\No{ \int_0^T U(\ta_1)^* V(\ta_1)Q(\ta_1) U(\ta_1) d\ta_1  }_{\FS^2} \\
	&\quad \ls \|V\|_{L^{4/3}_T L^2_x}^2 \|Q\|_{C_T \FS^2}.
\end{align}
Therefore, by Theorem \ref{th:OCK}, we get
\begin{align}
	A_3 \le C \|V\|_{L^{4/3}_T L^2_x}^2 \|Q\|_{C_T \FS^2} \le C T R^3 \le \frac{R}{8},
\end{align}
if we choose sufficiently small $T >0$.

\subsubsection{Estimate for $A_4$}
Since we have by Lemma \ref{lem:bilinear} and Kato--Seiler--Simon's inequality
\begin{align}
	&\No{\rh \Big(U(t) \int_0^T U(\ta)^* V(\ta)U(\ta) d\ta \int_0^T U(\ta_1)^* V(\ta_1)\ga_f U(\ta_1) d\ta_1 U(t)^*\Big)}_{L^4_T L^2_x} \\
	&\quad \ls \No{ \int_0^T U(\ta)^* V(\ta)U(\ta) d\ta \int_0^T U(\ta_1)^* V(\ta_1)\ga_f U(\ta_1) d\ta_1 }_{\FS^{4/3}} \\
	&\quad \le \No{ \int_0^T U(\ta)^* V(\ta)U(\ta) d\ta }_{\FS^4}
	\No{ \int_0^T U(\ta_1)^* V(\ta_1)\ga_f(\ta_1) U(\ta_1) d\ta_1  }_{\FS^2} \\
	&\quad \ls \|f\|_{L^2_\xi} \|V\|_{L^{4/3}_T L^2_x} \|V\|_{L^1_T L^2_x}.
\end{align}
Therefore, by Theorem \ref{th:OCK}, we get
\begin{align}
	A_4 \le C \|V\|_{L^{4/3}_T L^2_x} \|V\|_{L^1_T L^2_x} \le CT^{\frac{5}{4}} R^3 \le \frac{R}{8},
\end{align}
if we choose sufficiently small $T >0$.

\subsubsection{Estimate for $A_5$}
Since we have the explicit formula of $\CL_1(g)$ (see \cite[Proposition 1]{2014LS} and its proof), we have
\begin{align}
	A_5 &= \|\CL_1(g)\|_{L^4_T L^2_x}
	=\No{\CF \Dk{ \CL_1 (g)}(t,\xi)}_{L^4_T L^2_\xi} \\
	&\ls \No{
		\wh{w}(\xi) 
		\int_0^t 
		\sin((t-\ta)|\xi|^2) \check{f}(2(t-\ta)\xi)\wh{\rh_Q}(\ta,\xi) d\ta}_{L^4_T L^2_\xi} \\
	&\le \|\wh{w}\|_{L^\I_\xi} \|f\|_{L^1_\xi} T^{\ft}\int_0^T \|\rh_Q(\ta)\|_{L^2_x} d\ta \le \|\wh{w}\|_{L^\I_\xi} \|f\|_{L^1_\xi} T \|\rh_Q\|_{L^4_T L^2_x}.
\end{align}
Hence, we obtain
\begin{align}
	A_5 \le C T R \le \frac{R}{8},
\end{align}
if we choose sufficiently small $T>0$.

\subsubsection{Conclusion}
From the above, we obtain
\begin{align}
	\|\Ph_2[Q,g]\|_{L^4_T L^2_x}
	&\le A_1 + A_2 + \cdots + A_5 \le \frac{R}{2} + \frac{R}{8} \times 4 \le R.
\end{align}
Therefore, $\Ph:E(T,R) \to E(T,R)$ is well-defined. 
We can immediately prove $\Ph$ is a contraction mapping by multilinearizing the above estimates.
The uniqueness of the solution follows from the same argument. \qed

\subsection{Proof of Lemma \ref{lem:2D}}
We use the Duhamel formula \eqref{eq:Duhamel}.
Define
\begin{align}
	E(T,R):= \{(Q,g) \in C_T \FS^{2} \times L^2_T L^2_x : \|Q\|_{C_T \FS^{2}} \le R, 
	\|g\|_{L^2_T L^2_x} \le R\},
\end{align}
where $R:= 2(\|Q_0\|_{\FS^{2}} + \|\rh(U(t)_\st Q_0)\|_{L^2_T L^2_x})$ and $T>0$ will be chosen later.
We prove that $\Ph : E(T,R) \ni (Q,g) \to (\Ph_1[Q,g],\Ph_2[Q,g]) \in E(T,R)$ is a contraction mapping.

Let $(Q,g) \in E(T,R)$. Then we have by Lemma \ref{lem:bilinear}
\begin{align}
	\|\Ph_1[Q,g](t)\|_{\FS^{2}}
	&\le \|Q_0\|_{\FS^{2}}
	+ 2\No{\int_0^t U(\ta)^*V(\ta)Q(\ta)U(\ta) d\ta }_{\FS^{2}} \\
	&\quad + 2 \No{\int_0^t U(\ta)^*V(\ta)\ga_fU(\ta) d\ta }_{\FS^{2}} \\
	&\le \|Q_0\|_{\FS^{2}}
	+ C\|V\|_{L^\mu_T L^\nu_x} \|Q\|_{C_T \FS^{2}}
	+ C \|f\|_{L^1_\xi \cap L^\I_\xi} \|V\|_{L^1_T L^{2}_x} \\
	&\le \|Q_0\|_{\FS^{2}}
	+ \|w\|_{L^1_x \cap L^{1+\ep}_x} T^\de \|g\|_{L^{2}_T L^2_x} \|Q\|_{C_T \FS^{2}}
	+ C(f,w) T^\tw \|g\|_{L^2_T L^2_x} \\
	&\le \frac{R}{2} + CT^\de R^2 + CT^\tw R \le R,
\end{align}
for sufficiently small $T>0$, where we chose appropriate $\mu<2, \nu>2$ and $\de>0$.
For the density function, we have
\begin{align}
	\|\Ph_2[Q,g]\|_{L^2_T L^2_x}
	&\le \|\rh(U(t)_\st Q_0)\|_{L^2_T L^2_x}
	+ \|\rh(\CD_V[U(t)_\st Q_0])\|_{L^2_T L^2_x} \\
	&\quad + \|\rh(\CD_V^2[Q])\|_{L^2_T L^2_x}
	+ \|\rh(\CD_V^2[\ga_f])\|_{L^2_T L^2_x} + \|\CL_1(\rh_Q)\|_{L^2_T L^2_x} \\
	&=:A_1 + A_2 + A_3 + A_4+A_5.
\end{align}

\subsubsection{Estimate for $A_1$}
By the definition of $R$, we have
\begin{align}
	A_1 \le \frac{R}{2}.
\end{align}

\subsubsection{Estimate for $A_2$}
Let
\begin{align}
	\CV:= \int_0^T U(\ta)^* V(\ta) U(\ta) d\ta.
\end{align}
Note that
\begin{align}
	\|\rh(U(t)[\CV,Q_0]U(t)^*)\|_{L^2_T L^{2}_x}
	&\ls 
	\|\CV Q_0 \|_{\FS^{4/3}}\le
	\|\CV\|_{\FS^{4+}} \|Q_0\|_{\FS^{2-}} \ls \|V\|_{L^{\mu}_T L^{\nu}_x} \|Q_0\|_{\FS^{2-}}.
\end{align}
By the Christ--Kiselev lemma, we obtain
\begin{align}
	A_2 \le C T^\de \|w\|_{L^1_x \cap L^{1+\ep}_x} \|g\|_{L^{2}_T L^2_x} \|Q_0\|_{\FS^{2-0}}
	\le C T^\de R^2 \le \frac{R}{8}
\end{align}
for small $\de > 0$, if we choose sufficiently small $T > 0$.
When $Q_0 \in \CH^{0+}$, we can bound $A_2$ in the same way by combining Lemma \ref{lem:DP} and Theorem \ref{th:Ostri}.

\subsubsection{Estimate for $A_3$}
Note that
\begin{align}
	&\No{\rh \Big(U(t) \int_0^T U(\ta)^* V(\ta)U(\ta) d\ta \int_0^T U(\ta_1)^* V(\ta_1)Q(\ta_1) U(\ta_1) d\ta_1 U(t)^*\Big)}_{L^2_T L^2_x} \\
	&\quad \ls \No{ \int_0^T U(\ta)^* V(\ta)U(\ta) d\ta \int_0^T U(\ta_1)^* V(\ta_1)Q(\ta_1) U(\ta_1) d\ta_1 }_{\FS^{4/3}} \\
	&\quad \le \No{ \int_0^T U(\ta)^* V(\ta)U(\ta) d\ta }_{\FS^4}
	\No{ \int_0^T U(\ta_1)^* V(\ta_1)Q(\ta_1) U(\ta_1) d\ta_1  }_{\FS^2} \\
	&\quad \ls \|V\|_{L^2_T L^2_x} \|V\|_{L^\mu_T L^\nu_x} \|Q\|_{C_T \FS^2}.
\end{align}
Therefore, by Theorem \ref{th:OCK}, we get
\begin{align}
	A_3 \le C \|V\|_{L^2_T L^2_x} \|V\|_{L^\mu_T L^\nu_x} \|Q\|_{C_T \FS^2}
	\le C T^\de R^3 \le \frac{R}{8},
\end{align}

\subsubsection{Estimate for $A_4$}
We have
\begin{align}
	&\No{\rh \Big(U(t) \int_0^T U(\ta)^* V(\ta)U(\ta) d\ta \int_0^T U(\ta_1)^* V(\ta_1)\ga_f U(\ta_1) d\ta_1 U(t)^*\Big)}_{L^2_T L^2_x} \\
	&\quad \ls \No{ \int_0^T U(\ta)^* V(\ta)U(\ta) d\ta \int_0^T U(\ta_1)^* V(\ta_1)\ga_f U(\ta_1) d\ta_1 }_{\FS^{4/3}} \\
	&\quad \le \No{ \int_0^T U(\ta)^* V(\ta)U(\ta) d\ta }_{\FS^4}
	\No{ \int_0^T U(\ta_1)^* V(\ta_1)\ga_f(\ta_1) U(\ta_1) d\ta_1  }_{\FS^2} \\
	&\quad \ls \|f\|_{L^2_\xi} \|V\|_{L^{2}_T L^2_x} \|V\|_{L^1_T L^2_x}.
\end{align}
Therefore, by Theorem \ref{th:OCK}, we get
\begin{align}
	A_4 \le C T^{\tw} \|g\|_{L^{2}_T L^{2}_x}^2 \|Q\|_{C_T \FS^{2}} \le C T^\tw R^3 \le \frac{R}{8},
\end{align}
if we choose sufficiently small $T >0$.

\subsubsection{Estimate for $A_5$}
Since we have the explicit formula of $\CL_1(g)$ (see \cite[Proposition 1]{2014LS} and its proof), we have
\begin{align}
	A_5 &= \|\CL_1(g)\|_{L^2_T L^{2}_x}
	=\No{\CF [\CL_1 (g)](t,\xi)}_{L^2_T L^{2}_\xi} \\
	&\ls \No{
		\wh{w}(\xi) 
		\int_0^t 
		\sin((t-\ta)|\xi|^2) \check{f}(2(t-\ta)\xi)\wh{\rh_Q}(\ta,\xi) d\ta}_{L^2_T L^2_\xi} \\
	&\le \|\wh{w}\|_{L^2_\xi} \|\check{f}\|_{L^\I_\xi} T^{\tw}\int_0^T \|\rh_Q(\ta)\|_{L^2_x} d\ta \le C T \|\rh_Q\|_{L^2_T L^2_x} \le \frac{R}{8}.
\end{align}

\subsubsection{Conclusion}
From the above, we obtain
\begin{align}
	\|\Ph_2[Q,g]\|_{L^2_T L^2_x}
	&\le A_1 + A_2 + \cdots + A_5 \le \frac{R}{2} + \frac{R}{8} \times 4 \le R.
\end{align}
Therefore, $\Ph:E(T,R) \to E(T,R)$ is well-defined. 
We can immediately prove $\Ph$ is a contraction mapping by multilinearizing the above estimates.
The uniqueness of the solution follows from the same argument. 
\qed

\subsection{Proof of Lemma \ref{lem:3D}}
We use the Duhamel form:
\begin{align} \label{eq:Duhame2}
	\left \{
	\begin{aligned}
		Q(t) &= \Ph_1[Q,\rh_Q], \\
		V(t) &= w \ast \Ph_2[Q,V],
	\end{aligned}
	\right.
\end{align}
where
\begin{align}
	\Ph_1[Q,V]&:= U(t)_\st Q_0 + \CD_V[Q] +\CD_V[\ga_f]
\end{align}
and
\begin{align}
	\Ph_2[Q,V](t)
	&:= \rh(U(t)_\st Q_0) + \rh(\CD_V[U(t)_\st Q_0]) + \rh(\CD_V^2[U(t)_\st Q_0]) \\
	&\quad + \rh(\CD_V^3[Q]) + \rh(\CD_V^3[\ga_f]) + \rh(\CD_V^2[\ga_f]) + \rh(\CD_V[\ga_f]).
\end{align}
Define
\begin{align}
	E(T,R):= \{(Q,V) \in C_T \FS^{2} \times L^2_T (L^2_x \cap L^\I_x) : \|Q\|_{C_T \FS^{2}} \le R, 
	\|V\|_{L^2_T(L^2_x \cap L^\I_x)}\le R\},
\end{align}
where $R:= 2(\|Q_0\|_{\FS^{2}} + \|w \ast \rh(U(t)_\st Q_0)\|_{L^2_T(L^2_x \cap L^\I_x)})$ and $T>0$ will be chosen later.
We prove that $\Ph : E(T,R) \ni (Q,V) \to (\Ph_1[Q,V],w \ast \Ph_2[Q,V]) \in E(T,R)$ is a contraction mapping.

Let $(Q,V) \in E(T,R)$. Then we have by Kato--Seiler--Simon's inequality
\begin{align}
	\|\Ph_1[Q,V](t)\|_{\FS^{2}}
	&\le \|Q_0\|_{\FS^{2}}
	+ 2\No{\int_0^t U(\ta)^*V(\ta)Q(\ta)U(\ta) d\ta }_{\FS^{2}} \\
	&\quad + 2 \No{\int_0^t U(\ta)^*V(\ta)\ga_f(\ta)U(\ta) d\ta }_{\FS^{2}} \\
	&\le \|Q_0\|_{\FS^{2}}
	+ C\|V\|_{L^1_T L^\I_x} \|Q\|_{C_T \FS^{2}}
	+ C \|f\|_{L^1_\xi \cap L^\I_\xi} \|V\|_{L^1_T L^{2}_x} \\
	&\le \|Q_0\|_{\FS^{2}}
	+ CT^\tw \|V\|_{L^{2}_T(L^2_x \cap L^\I_x)} \|Q\|_{C_T \FS^{2}}
	+ C(f)T^\tw \|V\|_{L^2_T(L^2_x \cap L^\I_x)} \\
	&\le \frac{R}{2} + C T^\tw R^2 + C(f)T^\tw R \le R,
\end{align}
for sufficiently small $T>0$.
For the density function, we have
\begin{align}
	&\| w \ast \Ph_2[Q,g]\|_{L^2_T (L^2_x \cap L^\I_x)} \\
	&\quad \le \|w \ast \rh(U(t)_\st Q_0)\|_{L^2_T (L^2_x \cap L^\I_x)}
	+ C\|\rh(\CD_V[U(t)_\st Q_0])\|_{L^2_T L^{3/2}_x} \\
	&\qquad + C\|\rh(\CD_V^2[U(t)_\st Q_0])\|_{L^2_T L^{3/2}_x}
	+ C\|\rh(\CD_V^3[Q])\|_{L^2_T L^{3/2}_x} \\
	&\qquad + C\|\rh(\CD_V^3[\ga_f])\|_{L^2_T(L^1_x+L^2_x)}
	+ C\|\rh(\CD_V^2[\ga_f])\|_{L^2_T (L^1_x+L^2_x)} + C\|\rh(\CD_V[\ga_f])\|_{L^2_T L^{2}_x} \\
	&\quad =:A_1 + A_2 + A_3 + A_4+A_5+A_6+A_7.
\end{align}

\subsubsection{Estimate for $A_1$}
By the definition of $R$, we have
\begin{align}
	A_1 \le \frac{R}{2}.
\end{align}

\subsubsection{Estimate for $A_2$}
Define
\begin{align}
	\CV:= \int_0^T U(\ta)^* V(\ta) U(\ta) d\ta.
\end{align}
Note that
\begin{align}
	\|\rh(U(t)[\CV,Q_0]U(t)^*) \|_{L^2_T L^{3/2}_x}
	&\ls \|\CV Q_0 \|_{\FS^{6/5}}\le
	\| \CV \|_{\FS^{6+}} \|Q_0\|_{\FS^{3/2-}} \ls \|V\|_{L^{\mu}_T L^{\nu}_x} \|Q_0\|_{\FS^{3/2-}},
\end{align}
where $\mu<2$ and $\nu>3$.
By Christ--Kiselev's lemma, we obtain
\begin{align}
	A_2 \le C T^\de \|V\|_{L^{2}_T (L^{2}_x\cap L^\I_x)} \|Q_0\|_{\FS^{3/2-}}
	\le C\|Q_0\|_{\FS^{3/2-}} T^\de R \le \frac{R}{6}, 
\end{align}
for (small) $\de>0$, if we choose sufficiently small $T > 0$.
When $Q_0 \in \CH^{0+,\frac{3}{2}}$, we can bound $A_2$ in the same way by combining Lemma \ref{lem:DP} and Theorem \ref{th:Ostri}.

\subsubsection{Estimate for $A_3$}
Note that
\begin{align}
	&\|\rh (U(t) \CV^2 Q_0 U(t)^*)\|_{L^2_T L^{3/2}_x} \ls \| \CV^2 Q_0 \|_{\FS^{6/5}} \le \|\CV \|_{\FS^{12}}^2
	\|Q_0\|_{\FS^{3/2}}
	\ls \|V\|_{L^{4/3}_T L^{6}_x}^2 \|Q_0\|_{\FS^{3/2}}.
\end{align}
Hence, by Theorem \ref{th:OCK}, we have
\begin{align}
	A_3 \le C\|Q_0\|_{\FS^{3/2}} T^\tw R^2 \le \frac{R}{12}
\end{align}
for sufficiently small $T>0$.

\subsubsection{Estimate for $A_4$}
Next, we estimate $A_4$.
$A_4$ includes many terms, but we only consider a specific term; we can estimate the other terms in the same way.
Let $\phi \in C_c((0,T)\times \BR^3)$.
Let $\phi_0:= |\phi|^\tw$, $\phi =\phi_0 \phi_1$ and $V_0 := |V|^\tw$, $V=V_0 V_1$.
Then we have
\begin{align}
	&I := \bigg| \int_0^T dt \int_{\BR^d} dx \phi(x) \rh \bigg[U(t)\int_0^t d\ta U(\ta)^* V(\ta)U(\ta) d\ta \int_0^\ta d\ta_1 U(\ta_1)V(\ta_1)U(\ta_1)\\
	&\qquad \times \int_0^{T} d\ta_2 U(\ta_2)^* V(\ta_2) Q(\ta_2) U(\ta_2)  U(t)^* \bigg] \bigg|\\
	&\quad = \bigg|\Tr \bigg[ \int_0^T U(t)^* \phi(t) U(t) \int_0^t d\ta U(\ta)^* V(\ta) U(\ta) \\ 
	&\qquad \times \int_0^{\ta} d\ta_1 U(\ta_1) V(\ta_1) U(\ta_1) \int_0^T d\ta_2 U(\ta_2) V(\ta_2) Q(\ta_2) U(\ta_2)\bigg] \bigg| \\
	&\quad \le \|U^\oast \phi_0\|_{\FS^{12}_{\tx \to x}} \No{\phi_1(t)U(t) \int_0^t d\ta U(\ta)^* V_0(\ta)}_{\FS^6_{t,x}} \\
	&\qquad \times \No{V_1(\ta) U(\ta)\int_0^{\ta} U(\ta_1)^* V_0(\ta_1)}_{\FS^6_{t,x}}
	       \|V_1 U^\oc\|_{\FS^{12}_{x \to \tx}} \No{\int_0^T d\ta_2 U(\ta_2) V(\ta_2) Q(\ta_2) U(\ta_2) }_{\FS^2}.
\end{align}
By combining Lemma \ref{lem:DP} and Theorem \ref{th:Ostri}, we obtain
\begin{align}
	\|U^\oast \phi_0\|_{\FS^{12}_{\tx \to x}} \ls \|\ph\|_{L^2_t L^3_x}^\tw,
	\quad \|V_1 U^\oc\|_{\FS^{12}_{x \to \tx}} \ls \|V\|_{L^2_t L^3_x}^\tw.
\end{align}
In the same way, we have
\begin{align}
	\No{\phi_1(t)U(t) \int_0^T d\ta U(\ta)^* V_0(\ta)}_{\FS^6_{t,x}} \ls \|\phi\|_{L^2_t L^3_x}^\tw \|V\|_{L^2_t L^3_x}^\tw.
\end{align}
Hence, Theorem \ref{SCK} implies
\begin{align}
	\No{\phi_1(t)U(t) \int_0^t d\ta U(\ta)^* V_0(\ta)}_{\FS^6_{t,x}} \ls \|\phi\|_{L^2_t L^3_x}^\tw \|V\|_{L^2_t L^3_x}^\tw.
\end{align}
Collecting the above all, we conclude that
\begin{align}
	I \ls \|\phi\|_{L^2_t L^{3}_x} \|V\|_{L^2_t L^3_x}^2 \|V\|_{L^1_t L^\I_x} \|Q\|_{C_T \FS^2},
\end{align}
which implies
\begin{align}
	&\bigg\|\rh \bigg[U(t)\int_0^t d\ta U(\ta)^* V(\ta)U(\ta) d\ta \int_0^\ta d\ta_1 U(\ta_1)V(\ta_1)U(\ta_1) \\
	&\quad \times \int_0^{T} d\ta_2 U(\ta_2)^* V(\ta_2) Q(\ta_2) U(\ta_2)  U(t)^* \bigg] \bigg\|_{L^2_T L^{3/2}_x}
	\ls \|V\|_{L^2_t L^3_x}^2 \|V\|_{L^1_t L^\I_x} \|Q\|_{C_T \FS^2}.
\end{align}
By the standard argument to prove Theorem \ref{th:OCK} by the Whitney decomposition, we have
\begin{align}
	&\bigg\|\rh \bigg[U(t)\int_0^t d\ta U(\ta)^* V(\ta)U(\ta) d\ta \int_0^\ta d\ta_1 U(\ta_1)V(\ta_1)U(\ta_1) \\
	&\quad \times \int_0^{\ta_1} d\ta_2 U(\ta_2)^* V(\ta_2) Q(\ta_2) U(\ta_2)  U(t)^* \bigg] \bigg\|_{L^2_T L^{3/2}_x}\le  \|V\|_{L^2_t L^3_x}^2 \|V\|_{L^1_t L^\I_x} \|Q\|_{C_T \FS^2}.
\end{align}
By applying the same argument to the other terms, we get
\begin{align}
	A_4 &\le C T^\tw \|V\|_{L^2_T (L^2_x \cap L^\I_x)}^3 \|Q\|_{C_T \FS^2}
	\le CT^\tw R^4 \le \frac{R}{12},
\end{align}
for sufficiently small $T>0$. 

\subsubsection{Estimate for $A_5, A_6$}
Let $f_0:= |f|^\tw$ and $f=f_0f_1$.
Then we have
\begin{align}
	\|\rh(U(t) \CV^2 \ga_f U(t)^*)\|_{L^2_T (L^1_x + L^2_x)} 
	& \le\||\rh(U(t) \CV^2 \ga_{f_0^2} \CV^2 U(t)^*)|^\tw |\rh(\ga_{f_1^2})|^\tw \|_{L^2_T L^2_x} \\
	&\le C \|\rh(U(t) \CV^2 \ga_{f_0^2} \CV^2 U(t))\|_{L^1_T L^1_x}^\tw \\
	&\le C T^\tw \|\CV^2 \ga_{f_0} \ga_{f_0} \CV^2\|_{\FS^1}^\tw \\
	&\le C T^\tw \|\CV\|_{\CB(L^2_x)} \|\ga_{f_0} \CV\|_{\FS^2}\le C T^{\frac{3}{2}} \|V\|_{L^2_T(L^2_x \cap L^\I_x)}^2.
\end{align}
Moreover, we have
\begin{align}
	\|\rh(U(t) \CV \ga_f \CV U(t)^*)\|_{L^2_T(L^1_x + L^2_x)} 
	&\le T^\tw \|\CV \ga_{f_0} \ga_{f_1} \CV\|_{\FS^1} \\
	&\le C T^{\frac{3}{2}} \|V\|_{L^2_T(L^2_x \cap L^\I_x)}^2.
\end{align}
Therefore, by Theorem \ref{th:OCK}, we obtain
\begin{align}
	A_6 \le C(f) T^{\frac{3}{2}} \|V\|_{L^2_T(L^2_x \cap L^\I_x)}^2 \le \frac{R}{12}
\end{align}
for sufficiently small $T>0$.
Similarly, we have
\begin{align}
	A_5 \le \frac{R}{12}.
\end{align}

\subsubsection{Estimate for $A_7$}
Since we have the explicit formula of $\CL_1[\rh_Q]=-\rh(\CD_V[\ga_f])$ (see \cite[Proposition 1]{2014LS} and its proof), we have
\begin{align}
	A_4 &= \|\rh(\CD_V[\ga_f])\|_{L^2_T L^2_x} \\
	&\ls \No{
		\int_0^t 
		\sin((t-\ta)|\xi|^2) \check{f}(2(t-\ta)\xi)\wh{V}(\ta,\xi) d\ta}_{L^2_T L^2_\xi} \\
	&\le \|f\|_{L^1_\xi} T^{\tw}\int_0^T \|V(\ta)\|_{L^2_x} d\ta \\
	&\le C(f) T \|V\|_{L^2_T (L^2_x\cap L^\I_x)} \le \frac{R}{12}
\end{align}
for sufficiently small $T>0$.

\subsubsection{Conclusion}
From the above, we obtain
\begin{align}
	\|w \ast \Ph_2[Q,V]\|_{L^2_T(L^2_x \cap L^\I_x)}
	&\le A_1 + \cdots + A_7 \le \frac{R}{2} + \frac{R}{12}\times 6 = R.
\end{align}
Therefore, $\Ph:E(T,R) \to E(T,R)$ is well-defined. 
We can immediately prove $\Ph$ is a contraction mapping by multilinearizing the above estimates.
The uniqueness of the solution follows from the same argument. 
\qed

\section{Applications to the linearized Hartree equation: Scattering theory}
In this section, we prove Corollaries \ref{cor:singular scattering} and \ref{cor:scattering}.

\subsection{Proof of Corollary \ref{cor:singular scattering}} \label{sec:singular scattering proof}
Let $Q_0 \in \CH^\si$.
By Theorem \ref{prop:a.s.linear}, we have
\begin{align}
	\|\sd^\si \rh(U(t)_\st Q_0^{\om;\si})\|_{L^r_{\om} L^2_t L^2_x} \ls r^{\frac{1}{2}} \|Q_0\|_{\CH^\si}.
\end{align}
Therefore, for almost all $\om \in \Om$,
\begin{align}
	\|\rh(U(t)_\st Q_0^{\om;\si})\|_{L^2_t H^\si_x} <\I.
\end{align}
The Duhamel integral of \eqref{eq:singful random LQ} is 
\begin{align}\label{eq:L Duhamel}
	Q(t)= U(t)_\st Q_0^{\om;\si} - i \int_0^t U(t-\ta)_\st [w \ast \rh_Q(\ta),\ga_f] d\ta;
\end{align}
hence, we have
\begin{align}
	\rh_Q= \rh(U(t)_\st Q_0^{\om;\si}) - \CL_1(\rh_Q).
\end{align}
Therefore, we get
\begin{align}
	\rh_Q = (1+\CL_1)^{-1} \rh(U(t)_\st Q_0^{\om;\si}) \in L^2_t(\BR_+;H^\si_x)
\end{align}
because $\|(1+\CL_1)^{-1}\|_{\CB(L^2_t(\BR_+;L^2_x))} = \|(1+\CL_1)^{-1}\|_{\CB(L^2_t(\BR_+;H^\si_x))}$.
For this fact, see \cite[Proposition 1]{2014LS}.
Note that the uniqueness follows from the above argument.

We obtain the solution $Q(t)$ by \eqref{eq:L Duhamel} at least formally.
Now we justify it:
\begin{lem.}\label{lem:recover lemma}
	Let $d \ge 2$ and $\si = \frac{d}{2}-1$.
	Let $\lxr^{2\si} f(\xi) \in L^1_\xi \cap L^\I_\xi$, $Q_0 \in \CH^\si$ and $V \in L^2_t(\BR_+;H^\si_x)$.
	Then $Q(t)$ defined by \eqref{eq:L Duhamel} is in $C(\BR_+;\CH^\si)$.
	Furthermore, $Q(t)$ scatters; that is,
	there exist $Q_+ \in \FS^{\frac{2d}{d-1}}$ such that
	\begin{align}
		U(-t)Q(t)U(t) \to Q_+ \mbox{ in } \FS^{\frac{2d}{d-1}} \mbox{ as } t \to \I.
	\end{align}
\end{lem.}

\begin{proof}
First, we prove $Q(t) \in \CH^\si$ for all $t \in \BR_+$. 
By Kato--Seiler--Simon's inequality (\cite{1975SS}; see also \cite[Theorem 4.1]{2005Simon}) , we have
\begin{align}
	\|Q(t)\|_{\CH^\si}\ls \|Q_0^{\om;\si}\|_{\CH^{\si}} + \int_0^t \|V(\ta)\|_{H^\si_x} \|\lxr^{2\si} f(\xi)\|_{L^2_\xi} d\ta <\I.
\end{align}
It is easy to verify that $\BR \ni t \mapsto Q(t) \in \CH^\si$ is continuous.
Therefore, we get $Q(t) \in C(\BR_+;\CH^{\si})$.
Next, we prove the scattering.
Let $V_0:= |V|^{\frac{2}{d+2}}$ and $V=V_0 V_1$.
We have
\begin{align}
	&\|U(-t)Q(t)U(t)-U(-s)Q(s)U(s)\|_{\FS^{\frac{2d}{d-1}}} \\
	&\quad \ls \No{\int_s^t U(\ta)V(\ta)\ga_fU(\ta) d\ta}_{\FS^{\frac{2d}{d-1}}}  \\
	&\quad \le \No{\int_s^t U(\ta)V_0(\ta)}_{\FS^{2(d+2)}_{\tx\to x}} \|V_1(\ta) U(\ta) \ga_f\|_{\FS^{2d(d+2)/(d^2-2)}_{x \to \tx}} \\
	&\quad \ls \|V\|_{L^2_\ta(s,t;L^2_x)} \|\lxr^\si f(\xi)\|_{L^1_\xi \cap L^\I_\xi}
	\to 0 \mbox{ as } t,s \to \I.
\end{align}
To bound $\|V_1(\ta) U(\ta) \ga_f\|_{\FS^{2d(d+2)/(d^2-2)}_{x \to \tx}}$, we used Lemma \ref{lem:DP} and Theorem \ref{th:Ostri}.
By Theorem \ref{th:Ostri}, we have
\begin{align}
	\|\rh([U(t)\sd^{-\si}]_\st \ga_0)\|_{L^{(d+2)/2}_{t,x}} \ls \|\ga_0\|_{\FS^{d(d+2)/(2d+2)}}.
\end{align}
Therefore, by Lemma \ref{lem:DP}, we have
\begin{align}
	\|g U^\oc \sd^{-\si}\|_{\FS^{2d(d+2)/(d^2-2)}_{x \to \tx}} \ls \|g\|_{L^{2(d+2)/d}_{t,x}}.
\end{align}
From the above, we obtain
\begin{align}
	\|V_1(\ta)U(\ta)\ga_f\|_{{\FS^{2d(d+2)/(d^2-2)}_{x \to \tx}}}
	&\ls \|V_1 U^\oc \sd^{-\si}\|_{{\FS^{2d(d+2)/(d^2-2)}_{x \to \tx}}} \|\sd^\si \ga_f\|_{\CB(L^2_x)} \\
	&\le \|V_1\|_{L^{2(d+2)/d}_{t,x}} \|\lxr^\si f(\xi)\|_{L^1_\xi \cap L^\I_\xi}.
\end{align}
\end{proof}

\subsection{Proof of Corollary \ref{cor:scattering}}
For any $Q_0 \in \FS^2$, Theorem \ref{prop:a.s.linear2} implies
\begin{align}
	\|\rh(U(t)_\st Q_0^{\om,\wo})\|_{L^2_t L^2_x} <\I
\end{align}
for almost all $(\om, \wo) \in \Om\times \wt{\Om}$.
By the same argument as in Section \ref{sec:singular scattering proof}, we have
\begin{align}
	\rh_Q = (1+\CL_1)^{-1}\rh(U(t)_\st Q_0^{\om,\wo}) \in L^2_t(\BR_+;L^2_x)
\end{align}
for almost all $(\om, \wo) \in \Om\times \wt{\Om}$.
Hence, it suffices to show the following lemma, but we can prove it in the completely same way as the proof of Lemma \ref{lem:recover lemma}
\begin{lem.}
	Let $d \ge 2$ and $\si = \frac{d}{2}-1$.
	Let $\lxr^{\si} f(\xi) \in L^1_\xi \cap L^\I_\xi$, $Q_0 \in \CH^\si$ and $V \in L^2_t(\BR_+;H^\si_x)$.
	Then $Q(t)$ defined by \eqref{eq:L Duhamel} is in $C(\BR_+;\FS^2)$.
	Furthermore, $Q(t)$ scatters; that is,
	there exist $Q_+ \in \FS^{\frac{2d}{d-1}}$ such that
	\begin{align}
		U(-t)Q(t)U(t) \to Q_+ \mbox{ in } \FS^{\frac{2d}{d-1}} \mbox{ as } t \to \I.
	\end{align}
\end{lem.}

\section*{Acknowledgments}
The authors would like to express their deepest gratitude to their advisors K. Nakanishi and N. Kishimoto for their valuable time and support.
The first author was supported by JST, the establishment of university fellowships towards the creation of science technology innovation, Grant Number JPMJFS2123.

%%%%%%%%%%%%%%%%%%%%%%%%%%%%%%
%%%%%%%%%%%%%%%%%%%%%%%%%%%%%%
%%%%%%%%%%%%%%%%%%%%%%%%%%%%%%


\begin{thebibliography}{99}
	\bibitem{B84}
	Barab, J.E.: Nonexistence of asymptotically free solutions for a nonlinear Schr\"odinger equation.
	J. Math. Phys. \textbf{25} (1984), 3270--3273.
	
	\bibitem{BOP2015a}
	B\'enyi, \'A., Oh, T., Pocovnicu, O.:
	On the probabilistic Cauchy theory of the cubic nonlinear Schr\"odinger equation on Rd, $d\ge3$.
	Trans. Amer. Math. Soc. Ser. B \textbf{2} (2015), 1--50.

	\bibitem{BOP2015b}
	B\'enyi, \'A., Oh, T., Pocovnicu, O.:
	Wiener randomization on unbounded domains and an application to almost sure well-posedness of NLS.
	Excursions in Harmonic Analysis, in: Appl. Numer. Harmon. Anal.,
	Birkh\"auser/Springer, \textbf{4} 2015, 3--25.

	\bibitem{BOP2019}
	B\'enyi, \'A., Oh, T., Pocovnicu, O.:
	Higher order expansions for the probabilistic local Cauchy theory of the cubic nonlinear Schr\"odinger equation on $\R^3$.
	Trans. Amer. Math. Soc. Ser. B \textbf{6} (2019), 114--160.

	\bibitem{2019BHLNS}
	Bez, N., Hong, Y., Lee, S., Nakamura, S., Sawano, Y.:
	On the Strichartz estimates for orthonormal systems of initial data with regularity.
	Adv. Math. {\bf 354} (2019).
	
	\bibitem{Bez et al 2025 TLMS}
	Bez, N., Kinoshita, S., Shiraki, S.:
	Boundary Strichartz estimates and pointwise convergence for orthonormal systems
	Trans. London Math. Soc. 11 (2024), no. 1, Paper No. e70002, 43 pp.
	
	\bibitem{2020BLN}
	Bez, N., Lee, S., Nakamura, S.:
	Maximal estimates for the Schrödinger equation with orthonormal initial data.
	Selecta Math. (N.S.)
	\textbf{26} (2020), no.4.
	
	\bibitem{2021BLN}
	Bez, N., Lee, S., Nakamura, S.:
	Strichartz estimates for orthonormal families of initial data and weighted oscillatory integral estimates.
	Forum Math. Sigma {\bf 9} (2021), Paper No. e1.
	
	\bibitem{B94}
	Bourgain, J.:
	Periodic nonlinear Schrödinger equation and invariant measures.
	Comm. Math. Phys. \textbf{166} (1994).

	\bibitem{B96}
	Bourgain, J.:
	Invariant measures for the 2D-defocusing nonlinear Schr\"odinger equation.
	Comm. Math. Phys. \textbf{176} (1996), 421--445.
	
	\bibitem{1974BDF}
	Bove, A., Da Prato, G., Fano, G.:
	An existence proof for the Hartree-Fock time-dependent problem with bounded two-body interaction.
	Comm. Math. Phys.
	{\bf 37} (1974), 183--191.
	
	\bibitem{1976BDF}
	Bove, A., Da Prato, G., Fano, G.:
	On the Hartree-Fock time-dependent problem.
	Comm. Math. Phys.
	{\bf 49} (1976), 25--33.

	\bibitem{2019Br}
	Brereton, J.:
	Almost sure local well-posedness for the supercritical quintic NLS.
	Tunisian J. Math. \textbf{1} no.3 (2019), 427--453.
	
	\bibitem{BTT2013}
	Burq, N., Thomann, L., Tzvetkov, N.:
	Long time dynamics for the one dimensional non linear Schr\"odinger equation.
	Ann. Inst. Fourier (Grenoble) \textbf{63} no.6 (2013), 2137--2198.

	\bibitem{BT08a}
	Burq, N., Tzvetkov, N.:
	Random data Cauchy theory for supercritical wave equations I: Local theory.
	Invent. Math. \textbf{173} no.3 (2008), 449--475.
	
	\bibitem{BT08b}
	Burq, N., Tzvetkov, N.:
	Random data Cauchy theory for supercritical wave equations II: a local existence results.
	Invent. Math. \textbf{173} no.3 (2008), 477--496.
	
	\bibitem{2003BS}
	Birman, M.S., Solomyak, M.:
	Double operator integrals in a Hilbert space.
	Integral Equations Operator Theory {\bf 47} (2003), no.2, 131--168.
	
	\bibitem{1976C}
	Chadam, J.M.:
	The time-dependent Hartree-Fock equations with Coulomb two-body interaction.
	\newblock Comm. Math. Phys. 
	{\bf 46} (1976), 99--104.
	
	\bibitem{2001CK}
	Christ, M., Kiselev, A.:
	Maximal functions associated to filtrations.
	J. Funct. Anal. {\bf 179} (2001), no.2, 409--425.
	
	\bibitem{2017CHP}
	Chen, T., Hong, Y., Pavlovi\'{c}, N.: Global well-posedness of the NLS system for infinitely many fermions. Arch. Ration. Mech. Anal. {\bf 224} (2017), no.1, 91--123.
	
	\bibitem{2018CHP}
	Chen, T., Hong, Y., Pavlovi\'{c}, N.:
	On the scattering problem for infinitely many fermions in dimensions $d \geq 3$ positive temperature.
	Ann. Inst. H. Poincar\'{e} Anal. Non Lin\'{e}aire {\bf 35} (2018), no. 2, 393--416.

  \bibitem{DLM19}
  Dodson, B., J. L\"uhrmann, J., Mendelson, D.:
  Almost sure local well-posedness and scattering for the 4D cubic nonlinear Schrödinger equation.
  Adv. Math. \textbf{347} (2019), 619--676.
    
  \bibitem{DLM20}
  Dodson, B., L\"uhrmann, J., Mendelson, D.:
  Almost sure scattering for the 4d energy-critical defocusing nonlinear wave equation with radial data.
  Am. J. Math. \textbf{142} (2020), 475--504.

	\bibitem{2014FLLS} Frank, R.L., Lewin, M., Lieb, E.H., Seiringer, R.:
	Strichartz inequality for orthonormal functions.
	J. Eur. Math. Soc. {\bf 16} (2014), no. 7, 1507--1526.
	
	\bibitem{2017FS} Frank, R.L., Sabin, J.:
	Restriction theorems for orthonormal functions, Strichartz inequalities, and uniform Sobolev estimates.
	Amer. J. Math. {\bf 139} (2017), no. 6, 1649--1691.

\bibitem{Ghosh Swain 2024}		
Ghosh, S., Swain, J.:
On the Schatten exponent in orthonormal Strichartz estimate for the Dunkl operators.
Anal. Math. Phys. 14 (2024), no. 6, Paper No. 111, 15 pp.

	\bibitem{1970Go-Kre}
	Gohberg, I.C., Kre\u{\i}n, M.G.:
	Theory and applications of Volterra operators in Hilbert space.
	Transl. Math. Monogr., Vol. {\bf 24},
	American Mathematical Society, Providence, RI, (1970).
	
	\bibitem{2023Ha}
	Hadama, S.:
	Asymptotic stability of a wide class of steady states for the Hartree equation for random fields.
	arXiv:2303.02907
	
\bibitem{2023H}
	Hadama, S.:
	Asymptotic stability of a wide class of stationary solutions for the Hartree and Schrödinger equations for infinitely many particles.
Ann. H. Lebesgue \textbf{8} (2025), 181--218.

\bibitem{Hadama Hong 2025}
Hadama, S., Hong, Y.:
Global well-posedness of the nonlinear Hartree equation for infinitely many particles with singular interaction.
J. Funct. Anal. \textbf{289} (2025), no. 9, Paper No. 111102, 38 pp.
	
	\bibitem{Hoshiya 2024 JFA}
	Hoshiya, A.:
	Orthonormal Strichartz estimate for dispersive equations with potentials.
	J. Funct. Anal. \textbf{286} (2024), no. 11, Paper No. 110425, 63 pp.
	
	\bibitem{Hoshiya 2025 JMP}
	Hoshiya, A.:
	Orthonormal Strichartz estimates for Schrödinger operator and their applications to infinitely many particle systems.(English summary)
	J. Math. Phys. \textbf{66} (2025), no. 7, Paper No. 071505, 19 pp.
	
	\bibitem{Hoshiya 2025 AHP}
	Hoshiya, A.:
	Uniform Resolvent and Orthonormal Strichartz Estimates for Repulsive Hamiltonian.
	Ann. Henri Poincar\'{e} (2025).
	
	\bibitem{1998KT}
	Keel, M., Tao, T.:
	Endpoint Strichartz estimates.
	Amer. J. Math. {\bf 120} (1998), no.5, 955--980.
	
	
	\bibitem{2015LS}
	Lewin, M., Sabin, J.: The Hartree equation for infinitely many particles I. Well-posedness theory. Comm. Math. Phys. {\bf 334} (2015), no.1, 117--170.
	
	
	\bibitem{2014LS}
	Lewin, M., Sabin, J.: The Hartree equation for infinitely many particles II: Dispersion and scattering in 2D. Anal. PDE {\bf 7} (2014), no. 6, 1339--1363.
	
	\bibitem{2014LM}
	L\"{u}hrmann, J., Mendelson, D.:
	Random data Cauchy theory for nonlinear wave equations of power-type on $\BR^3$.
	Comm. Partial Differential Equations \textbf{39} (2014), no.12, 2262--2283.
	
	\bibitem{Mondal Song 2025}			
	Mondal, S.S., Song, M.:
		Orthonormal Strichartz inequalities for the $(k,a)$-generalized Laguerre operator and Dunkl operator.
		Isr. J. Math. \textbf{269}, 697–729 (2025).
		
	\bibitem{M19}
	Murphy, J.:
	Random data final-state problem for the mass-subcritical NLS in $L^{2}$.
	Proc. Amer. Math. Soc. \textbf{147} no. 1 (2019), 339--350.
	
	\bibitem{2020N} Nakamura, S.: The orthonormal Strichartz inequality on torus. Trans. Amer. Math. Soc.
	{\bf 373} (2020), no. 2, 1455--1476.
	
	\bibitem{NY19}
	Nakanishi, K., Yamamoto, T.:
	Randomized Final-data Problem for Systems of Nonlinear Schr\"odinger Equations and the Gross-Pitaevskii Equation.
	Math. Res. Lett. \textbf{26} no. 1 (2019), 253--279.
	
	
	\bibitem{Nguyen You 2025}
	Nguyen, T.T., You, C.:
		Plasmons for the Hartree equations with Coulomb interaction.
		Probab. Math. Phys. \textbf{6} (2025), no. 3, 913--960.
	
	\bibitem{2016OP}
	Oh, T., Pocovnicu, O.:
	Probabilistic global well-posedness of the energy-critical defocusing quintic nonlinear wave equation on $\R^3$.
	J. Math. Pures Appl. \textbf{105} 3 (2016) 342--366.
	
	\bibitem{2021PS}
	Pusateri, F., Sigal, I.M.:
	Long-time behaviour of time-dependent density functional theory.
	Arch. Ration. Mech. Anal. {\bf 241} (2021), no. 1, 447--473.
	
	\bibitem{1975SS}
	Seiler, E., Simon, B.: Bounds in the Yukawa2 quantum field theory: upper bound on the pressure, Hamiltonian bound and linear lower bound.
	Comm. Math. Phys. {\bf 45}, 99--114 (1975).
	
	\bibitem{Senapati et al 2024}
	Senapati, P.J.K., Boggarapu, P., Mondal, S.S., Mejjaoli, H.:
		Restriction theorem for the Fourier-Dunkl transform and its applications to Strichartz inequalities.
		J. Geom. Anal. \textbf{34} (2024), no. 3, Paper No. 74, 36 pp.
	
	\bibitem{2005Simon}
	Simon, B.:
	{\it Trace ideals and their applications, Second edition.}
	Mathematical Surveys and Monographs. 
	{\bf 120}. American Mathematical Society. (2005)
	
	\bibitem{Sp22}
	Spitz, M.:
	Randomized final-state problem for the Zakharov system in dimension three.
	Comm. Partial Differential Equations \textbf{47} no.2 (2022) 346--377.

	\bibitem{2009Th}
	Thomann, L.:
	Random data Cauchy problem for supercritical Schr\"odinger equations.
	Ann. Inst. H. Poincaré Anal. Non Linéaire \textbf{26} no.6 (2009), 2385--2402.

	\bibitem{2008Tz}
	Tzvetkov, N.:
	Invariant measures for the defocusing nonlinear Schrödinger equation.
	Annales de l'Institut Fourier \textbf{58} no.7 (2008), 2543--2604.

	\bibitem{W32}
	Wiener, N.:
	Tauberian theorems.
	Ann. of Math. (2) \textbf{33} no.1 (1932), 1--100.

	\bibitem{1992Z}
	Zagatti, S.:
	The Cauchy problem for Hartree-Fock time-dependent equations.
	Ann. Inst. H. Poincar\'{e} Phys. Théor. {\bf 56} (1992), 357--374.
	
	\bibitem{ZF12}
	Zhang, T., Fang, D.:
	Random data Cauchy theory for the generalized incompressible Navier-Stokes equations.
	J. Math. Fluid Mech. \textbf{14} (2012), 311--324.
\end{thebibliography}
\end{document}